\documentclass[11pt]{article}
\RequirePackage[OT1]{fontenc}

\RequirePackage[aos,amsthm,amsmath,natbib,doublespacing]{imsart}
\usepackage{amsfonts,graphicx}
\newcommand {\apgt} {\ {\raise-.5ex\hbox{$\buildrel>\over\sim$}}\ }
\newcommand {\aplt} {\ {\raise-.5ex\hbox{$\buildrel<\over\sim$}}\ }
\begin{document}

\begin{frontmatter}

\title{The Penalized Profile Sampler}
\runtitle{The Penalized Profile Sampler}
\author{Guang Cheng\thanksref{t1}}
\and
\author{Michael R. Kosorok\thanksref{t1}}
\affiliation{Duke University and University of North Carolina at
Chapel Hill}
\address{Guang Cheng\\ Institute of Statistics and\\ Decision Sciences\\
Duke University\\
214 Old Chemistry Building\\
Durham, NC  27708\\USA\\Email: chengg@stat.duke.edu}

\address{Michael R. Kosorok\\ Department of Biostatistics\\
School of Public Health\\
University of North Carolina\\ at Chapel Hill\\
3101 McGavran-Greenberg Hall\\
 Chapel Hill, NC  27599\\USA\\Email: kosorok@unc.edu}
\runauthor{G. Cheng and M. R. Kosorok}
\thankstext{t1}{Supported in part by CA075142}

\begin{abstract}
The penalized profile sampler for semiparametric inference
is an extension of the profile
sampler method \cite{lkf05} obtained by profiling a penalized
log-likelihood. The idea is to base inference on the
posterior distribution obtained by multiplying a profiled
penalized log-likelihood by a prior for the parametric component,
where the profiling and penalization are applied to the nuisance
parameter. Because the prior is not applied to the full
likelihood, the method is not strictly Bayesian.  A benefit of
this approximately Bayesian method is that it circumvents the
need to put a prior on the possibly infinite-dimensional
nuisance components of the model.  We investigate the first and
second order frequentist performance of the penalized profile
sampler, and demonstrate that the accuracy of the procedure can be
adjusted by the size of the assigned smoothing parameter.
The theoretical validity of the procedure is illustrated for two
examples: a partly linear model with normal error for current
status data and a semiparametric logistic regression model.
As far as we are aware, there are no other methods of inference
in this context known to have second order frequentist validity.
\end{abstract}



\end{frontmatter}

\newtheorem{theorem}{\indent \sc Theorem}
\newtheorem{corollary}{\indent \sc Corollary}
\newtheorem{lemma}{\indent \sc Lemma}
\newtheorem{proposition}{\indent \sc Proposition}
\newtheorem{remark}{\indent \sc Remark}
\newcommand{\phif}{\textsc{igf}}
\newcommand{\sign}{\mbox{sign}}
\newcommand{\phgf}{\textsc{gf}}
\newcommand{\fix}{$\textsc{gf}_0$}
\newcommand{\mb}[1]{\mbox{\bf #1}}
\newcommand{\Exp}[1]{\mbox{E}\left[#1\right]}
\newcommand{\pr}[1]{\mbox{P}\left[#1\right]}
\newcommand{\pp}[0]{\mathbb{P}}
\newcommand{\ee}[0]{\mbox{E}}
\newcommand{\re}[0]{\mathbb{R}}
\newcommand{\argmax}[0]{\mbox{argmax}}
\newcommand{\argmin}[0]{\mbox{argmin}}
\newcommand{\ind}[0]{\mbox{\Large\bf 1}}
\newcommand{\narrow}{\stackrel{n\rightarrow\infty}{\longrightarrow}}
\newcommand{\weakpn}{\stackrel{P_n}{\leadsto}}
\newcommand{\weakpnboot}{\mbox{\raisebox{-1.5ex}{$\stackrel
{\mbox{\scriptsize $P_n$}}{\stackrel{\mbox{\normalsize
$\leadsto$}} {\mbox{\normalsize $\circ$}}}$}}\,}
\newcommand{\ol}[1]{\overline{#1}}
\newcommand{\avgse}[1] { \bar{\hat{\sigma}}_{#1} }
\newcommand{\mcse}[1]  { \sigma^{*}_{#1} }
\newcommand{\po}{\textsc{po}}
\newcommand{\ph}{\textsc{ph}}

\section{Introduction}
Semiparametric models are statistical models indexed by both a
finite dimensional parameter of interest $\theta$ and an infinite
dimensional nuisance parameter $\eta$. The profile likelihood is
typically defined as
\begin{eqnarray*}
pl_{n}(\theta)=\sup_{\eta\in\mathcal{H}}lik_{n}(\theta,\eta),
\end{eqnarray*}
where $lik_{n}(\theta,\eta)$ is the likelihood of the
semiparametric model given $n$ observations and $\mathcal{H}$ is
the parameter space for $\eta$. We also define
\begin{eqnarray*}
\hat{\eta}_{\theta}=\argmax_{\eta\in\mathcal{H}}lik_{n}(\theta,\eta).
\end{eqnarray*}
The convergence rate of the nuisance parameter $\eta$ is
the order of $d(\hat{\eta}_{\tilde{\theta}_{n}},\eta_{0})$, where
$d(\cdot,\cdot)$ is some metric on $\eta$,
$\tilde{\theta}_n$ is any sequence satisfying $\tilde{\theta}_n
=\theta_0+o_P(1)$, and $\eta_{0}$ is the true
value of $\eta$. Typically,
\begin{eqnarray}
d(\hat{\eta}_{\tilde{\theta}_{n}},\eta_{0})=O_{P}(\|\tilde{\theta}_{n}-\theta_{0}\|+n^{-r}),\label{convr}
\end{eqnarray}
where $\|\cdot\|$ is the Euclidean norm and $r>1/4$. Of course, a
smaller value of $r$ leads to a slower convergence rate of the
nuisance parameter. For instance, the nuisance parameter in the Cox
proportional hazards model with right censored data, the
cumulative hazard function, has the parametric rate, i.e., $r=1/2$. If
current status data is applied to the Cox model instead, then the
convergence rate will be slower, with $r=1/3$, due to the loss of
information provided by this kind of data.

The profile sampler is the procedure of sampling from the
posterior of the profile likelihood in order to estimate and
draw inference on the parametric component $\theta$ in a
semiparametric model, where the profiling is done over the
possibly infinite-dimensional nuisance parameter $\eta$.
\cite{lkf05} show that the profile sampler gives a
first order correct approximation to the maximum likelihood
estimator $\hat{\theta}_{n}$ and consistent estimation of the
efficient Fisher information for $\theta$ even when the nuisance
parameter is not estimable at the $\sqrt{n}$ rate. Another
Bayesian procedure employed to do semiparametric estimation
is considered in \cite{s02} who study the marginal semiparametric
posterior distribution for a parameter of interest. In particular,
\cite{s02} show that marginal semiparametric posterior
distributions are asymptotically normal and centered at the
corresponding maximum likelihood estimates or posterior means,
with covariance matrix equal to the inverse of the Fisher information.
Unfortunately, this fully Bayesian method requires specification of
a prior on $\eta$, which is quite challenging since for some
models there is no direct extension of the concept of a Lebesgue
dominating measure for the infinite-dimensional parameter set
involved \cite{k75}. The advantages of the profile sampler
for estimating $\theta$ compared to other methods is discussed
extensively in \cite{ck06}, \cite{ck062} and \cite{lkf05}.

In many semiparametric models involving a smooth nuisance parameter,
it is often convenient and beneficial to perform estimation using
penalization. One motivation for this is that, in the absence of any
restrictions on the form of the function $\eta$, maximum
likelihood estimation for some semiparametric models leads to
over-fitting. Seminal applications of penalized maximum likelihood estimation include estimation of a probability density function in \cite{si82} and nonparametric linear regression
in \cite{si85}. Note that penalized likelihood is a special case
of penalized quasi-likelihood studied in \cite{ma97}. Under certain
reasonable regularity conditions,
penalized semiparametric log-likelihood estimation
can yield fully efficient estimates for
$\theta$ (see, for example, \cite{ma97}). As far as we are aware,
the only general procedure for inference for $\theta$
in this context known to be theoretically valid is a weighted
bootstrap with bounded random weights (see \cite{mk05}).  It is even
unclear whether the usual nonparametric bootstrap will work in
this context when the nuisance parameter has a convergence
rate $r<1/2$.

In contrast,
\cite{ck06} and \cite{ck062} have shown that the profile sampler
procedure without penalization
can essentially yield second order frequentist valid
inference for $\theta$ in semiparametric models, where the estimation
accuracy is dependent on the convergence rate of the nuisance parameter. In other words, a faster convergence rate
of the nuisance parameters can yield more precise frequentist
inference for $\theta$. These second order results are
verified in \cite{ck06} and \cite{ck062} for several examples,
including the Cox model for both right censored and current
status data, the proportional odds model, case-control
studies with missing covariates, and the partly linear normal model.
The convergence rates for these
models range from the parametric to the cubic. The
work in \cite{ck062} has shown clearly that
the accuracy of the inference for $\theta$ based on the
profile sampler method is intrinsically determined by the
semiparametric model specifications through its entropy number.

The purpose of this paper is to ask the
somewhat natural question: does sampling from
a profiled penalized log-likelihood (which process we
refer hereafter to as the penalized profile sampler)
yield first and even
second order accurate frequentist inference? The conclusion
of this paper is that the answer is yes and, moreover,
the accuracy of the inference depends in a fairly simple way
on the size of the smoothing parameter.

The unknown parameters in the semiparametric models we study in
this paper includes $\theta$, which we assume belongs to some
compact set $\Theta\subset \mathbb{R}^{d}$, and $\eta$, which we
assume to be a function in the Sobolev class of functions supported
on some compact set on the real line, whose $k$-th derivative
exists and is absolutely continuous with $J(\eta)<\infty$, where
\begin{eqnarray*}
J^{2}(\eta)=\int_{\mathcal{Z}}(\eta^{(k)}(z))^{2}dz.
\end{eqnarray*}
Here $k$ is a fixed, positive integer and $\eta^{(j)}$ is the $j$-th
derivative of $\eta$ with respect to $z$.
Obviously $J^{2}(\eta)$ is some
measurement of complexity of $\eta$. We denote
$\mathcal{H}_{k}$ as the Sobolev function class with degree $k$.
The penalized log-likelihood in this context is:
\begin{eqnarray}
\log
lik_{\lambda_{n}}(\theta,\eta)=\log lik(\theta,\eta)
-\lambda_{n}^{2}J^{2}(\eta),\label{penlik}
\end{eqnarray}
where $\log lik(\theta,\eta)\equiv\mathbb{P}_{n}
\ell_{\theta,\eta}(X)$, $\ell_{\theta,\eta}(X)$ is the
log-likelihood of the single observation $X$, and
$\lambda_{n}$ is a smoothing parameter, possibly dependent on
data. In practice, $\lambda_{n}$ can be obtained by
cross-validation \cite{w90} or by inspecting the various curves
for different values of $\lambda_{n}$. The penalized maximum
likelihood estimators $\hat{\theta}_{n}$ and $\hat{\eta}_{n}$
depend on the choice of the smoothing parameter $\lambda_{n}$.
Consequently we use the notation $\hat{\theta}_{\lambda_{n}}$ and
$\hat{\eta}_{\lambda_{n}}$ for the remainder of this paper
to denote the estimators obtained from maximizing~(\ref{penlik}). In
particular, a larger smoothing parameter usually leads to a less
rough penalized estimator of $\eta_{0}$.

For the purpose of establishing first order accuracy of
inference for $\theta$ based on the penalized profile sampler,
we assume that the bounds for the smoothing parameter
are in the form below:
\begin{eqnarray}
\lambda_{n}=o_{P}(n^{-1/4}) \;\;\mbox{and}\;\;
\lambda_{n}^{-1}=O_{P}(n^{k/(2k+1)}).\label{smooth}
\end{eqnarray}
The condition (\ref{smooth}) is assumed to hold throughout this
paper. One way to ensure (\ref{smooth}) in practice is simply to
set $\lambda_{n}=n^{-k/(2k+1)}$. Or we can just choose
$\lambda_{n}=n^{-1/3}$ which is independent of $k$. It turns
out that the upper
bound guarantees that $\hat{\theta}_{\lambda_{n}}$ is
$\sqrt{n}$-consistent, while
the lower bound controls the penalized nuisance
parameter estimator convergence rate.
Another approach to controlling
estimators is to use sieve estimates with assumptions on the
derivatives (see \cite{h99}). We will not pursue this further here.

The log-profile penalized likelihood is defined as follows:
\begin{eqnarray}
\log pl_{\lambda_{n}}(\theta)=\log
lik(\theta,\hat{\eta}_{\theta,\lambda_{n}})-\lambda_{n}^{2}
J^{2}(\hat{\eta}_{\theta,\lambda_{n}})\label{propen},
\end{eqnarray}
where $\hat{\eta}_{\theta,\lambda_{n}}$ is
$argmax_{\eta\in\mathcal{H}_{k}}\log
lik_{\lambda_{n}}(\theta,\eta)$ for fixed $\theta$ and
$\lambda_{n}$. The penalized profile sampler is just the procedure
of sampling from the posterior distribution of
$pl_{\lambda_{n}}(\theta)$ by assigning a prior on $\theta$. By
analyzing the corresponding MCMC chain from the frequentist's
point of view, our paper obtains the following conclusions:
\begin{itemize}
\item[1] {\it Distribution Approximation}: The posterior
distribution with respect to $pl_{\lambda_{n}}(\theta)$ can be
approximated by the normal distribution with mean the maximum
penalized likelihood estimator of $\theta$ and variance the
inverse of the efficient information matrix, with error
$O_{P}(n^{1/2}\lambda_{n}^{2})$;
\item[2] {\it Moment Approximation}:
The maximum penalized likelihood estimator of
$\theta$ can be approximated by the mean of the MCMC chain with error
$O_{P}(\lambda_{n}^{2})$. The efficient information matrix can be
approximated by the inverse of the variance of the MCMC chain with
error $O_{P}(n^{1/2}\lambda_{n}^{2})$;
\item[3] {\it Confidence Interval Approximation}: An exact
frequentist confidence interval of Wald's type for $\theta$ can be
estimated by the credible set obtained from the MCMC chain
with error $O_{P}(\lambda_{n}^{2})$.
\end{itemize}

Obviously, given any smoothing parameter satisfying the upper
bound in~(\ref{smooth}), the penalized profile sampler can yield
first order frequentist valid inference for $\theta$, similar as
to what was shown for the profile sampler in \cite{lkf05}.
Moreover, the above conclusions are actually second order
frequentist valid results, whose approximation accuracy is directly
controlled by the smoothing parameter. Note that the corresponding
results for the usual (non-penalized) profile sampler with
nuisance parameter convergence rate $r$ in
\cite{ck062} are obtained by replacing in the above
$O_{P}(n^{1/2}\lambda_{n}^{2})$ with $O_{P}(n^{-1/2}\vee
n^{-r+1/2})$ and $O_{P}(\lambda_{n}^{2})$ with $O_{P}(n^{-1}\vee
n^{-r})$, for all respective occur
where $r$ is as defined in~(\ref{convr}).

Our results are the first higher order frequentist inference
results for penalized semiparametric estimation. The layout of
the article is as follows. The next section, section 2, introduces
the two main examples we will be using for illustration:
partly linear regression for current status data and
semiparametric logistic regression. Some background
is given in section~3, including
the concept of a least favorable submodel as well as
some notations and the main model assumptions.
In section~\ref{rate}, some preliminary results are developed,
including three rather different theorems
concerning the convergence rates of the penalized nuisance
parameters and the order of the estimated penalty term under
different conditions. The corresponding rates for the
two featured examples are also calculated in this section.
The main results and implications are discussed in section 5, and
all remaining model assumptions are verified for the
examples in section 6. A brief discussion of future work is
given in section~7. We postpone all technical tools and
proofs to the last section, section 8.

\section{Examples}

\subsection{Partly Linear Normal Model with Current Status Data}
In this example, we study the partly linear regression model with
normal residue error. The continuous outcome $Y$, conditional on
the covariates $(U,V)\in \mathbb{R}^{d}\times\mathbb{R}$, is
modeled as
\begin{eqnarray}
Y=\theta^{T}U+f(V)+\epsilon,\label{eg2lik}
\end{eqnarray}
where $f$ is an unknown smooth function, and $\epsilon\sim
N(0,\sigma^{2})$ with finite variance $\sigma^{2}$. For
simplicity, we assume for the rest of the paper that
$\sigma=1$. The theory we propose also
works when $\sigma$ is unknown, but the added complexity would
detract from the main issues.  We also assume that only the current
status of response $Y$ is observed at a random
censoring time $C\in\mathbb{R}$. In other words, we observe
$X=(C,\Delta,U,V)$, where indicator $\Delta=1\{Y\leq C\}$. Current
status data may occur due to study design or measurement
limitations. Examples of such data arise in several fields,
including demography, epidemiology and econometrics. For
simplicity of exposition, $\theta$ is assumed to be one
dimensional.

Under the model~(\ref{eg2lik}) and given that the joint distribution
for $(C,U,V)$ does not involve parameters $(\theta,f)$, the
log-likelihood for a single observation at
$X=x\equiv(c,\delta,u,v)$ is
\begin{eqnarray}
loglik_{\theta,f}(x)&=&\delta\log\left\{\Phi\left(c-\theta
u-f(v)\right)\right\}\nonumber\\&&+(1-\delta)\log\left\{1-\Phi\left(c-\theta
u-f(v)\right)\right\},\label{eg2p}
\end{eqnarray}
where $\Phi$ is the standard normal distribution. The parameter of
interest, $\theta$, is assumed to belong to some compact set in
$\mathbb{R}^{1}$. The nuisance parameter is the function $f$,
which belongs to the Sobolev function class of degree $k$. We
further make the following assumptions on this model. We assume
that $(Y,C)$ is independent given $(U,V)$. The covariates $(U,V)$
are assumed to belong to some compact set, and the support for
random censoring time $C$ is an interval $[l_{c},u_{c}]$, where
$-\infty<l_{c}<u_{c}<\infty$. In addition, $E Var(U|V)$ is strictly
positive and $Ef(V)=0$. The first order asymptotic behaviors of
the penalized log-likelihood estimates of a slightly more
general version of this model have been
extensively studied in \cite{mm05}.

\subsection{Semiparametric Logistic Regression}
Let $X_{1}=(Y_{1},W_{1},Z_{1})$, $X_{2}=(Y_{2},W_{2},Z_{2})$,
\ldots be independent copies of $X=(Y,W,Z)$, where $Y$ is a
dichotomous variable with conditional expectation
$E(Y|W,Z)=F(\theta^{T}
W+\eta(Z))$. $F(u)$ is the logistic distribution defined as
$e^{u}/(e^{u}+1)$. Obviously the likelihood for a single
observation is of the following form:
\begin{eqnarray}
p_{\theta,\eta}(x)=F(\theta^{T} w+\eta(z))^{y}(1-F(\theta^{T}
w+\eta(z)))^{1-y}f^{(W,Z)}(w,z).\label{eg1den}
\end{eqnarray}

This example is a special case of quasi-likelihood in partly
linear models when the conditional variance of response $Y$ is
taken to have some quadratic form of the conditional mean of $Y$.
In the absence of any restrictions on the form of the function
$\eta$, the maximum likelihood of this simple model often leads to
over-fitting. Hence \cite{gg71} propose maximizing instead the
penalized likelihood of the form $\log
lik(\theta,\eta)-\lambda_{n}^{2}J^{2}(\eta)$; and \cite{ma97}
studied the asymptotic properties of the maximum penalized
likelihood estimators for $\theta$ and $\eta$. For simplicity, we
will restrict ourselves to the case where
$\Theta\subset\mathbb{R}^{1}$ and $(W,Z)$ have bounded support,
say $[0,1]^{2}$. To ensure the identifiability of the parameters,
we assume that $E Var(W|Z)$ is positive and that the support of
$\mathcal{Z}$ contains at least $k$ distinct points in $[0,1]$.

\begin{remark}
Another interesting
potential example we may apply the penalized profile
sampler method to is the classic proportional hazards model with
current status data by penalizing the cumulative hazard function
with its Sobolev norm. There are two motivations for us to
penalize the cumulative hazard function in the Cox model. One is
that the estimated step functions from the unpenalized estimation
cannot be used easily for other
estimation or inference purposes. Another issue with the unpenalized approach is that without
making stronger continuity assumptions, we cannot achieve
uniform consistency even on a compact set \cite{mm05}. The
asymptotic properties of the corresponding penalized M-estimators
have been studied in \cite{mm06}.
\end{remark}

\section{Preliminaries}
In this section, we present some necessary preliminary material
concerning least favorable submodels, general notational conventions
for the paper, and an enumeration of the main assumptions.
\subsection{Least favorable submodels}
In this subsection, we briefly review the concept of a least
favorable submodel. A submodel $t\mapsto p_{t,\eta_{t}}$ is
defined to be {least favorable} at $(\theta,\eta)$ if
$\tilde{\ell}_{\theta,\eta}=
\partial/\partial t\log p_{t,\eta_{t}}$, given $t=\theta$, where
$\tilde{\ell}_{\theta,\eta}$ is the {efficient score function}
for $\theta$. The {efficient score function} for $\theta$ can
be viewed as the projection of the {score function} for
$\theta$ onto the {tangent space} of $\eta$. The inverse of
its variance is exactly the {efficient information matrix}
$\tilde{I}_{\theta,\eta}$. We abbreviate hereafter
$\tilde{\ell}_{\theta_{0},\eta_{0}}$ and
$\tilde{I}_{\theta_{0},\eta_{0}}$ with $\tilde{\ell}_{0}$ and
$\tilde{I}_{0}$, respectively. The ``direction'' along which
$\eta_{t}$ approaches $\eta$ in the least favorable submodel is
called the {least favorable direction}. An insightful review
about least favorable submodels and efficient score functions can
be found in Chapter~3 of \cite{kta}. By the above construction of
the least favorable submodel, $\log pl_{\lambda_{n}}(\theta)$ can
be rewritten in the following form:
\begin{eqnarray}
\log
pl_{\lambda_{n}}(\theta)&=&\ell(\theta,\theta,\hat{\eta}_{\theta,\lambda_{n}})-\lambda_{n}^{2}J^{2}(\eta_{\theta}(\theta,\hat{\eta}_{\theta,
\lambda_{n}})),\label{penlfs}
\end{eqnarray}
where $\ell(t,\theta,\eta)(x)=\log
lik(t,\eta_{t}(\theta,\eta))(x)$, $t\mapsto\eta_{t}(\theta, \eta)$
is a general map from the neighborhood of $\theta$ into the
parameter set for $\eta$, with $\eta_{\theta}(\theta,\eta)=\eta$.
The concrete forms of (\ref{penlfs}) will depend on the
situation.

\subsection{Notation}
We present in this subsection
some notation that will be used throughout the paper.
The derivatives of the function
$\ell(t,\theta,\eta)$ are with respect to its first argument, $t$.
For the derivatives relative to the other two arguments $\theta$
and $\eta$, we use the following shortened notation:
$\ell_{\theta}(t,\theta,\eta)$ indicates the first derivative of
$\ell(t,\theta,\eta)$ with respect to $\theta$. Similarly,
$\ell_{t,\theta}(t,\theta,\eta)$ denotes the derivative of
$\dot{\ell}(t,\theta,\eta)$ with respect to $\theta$. Also,
$\ell_{t,t}(\theta)$ and $\ell_{t,\theta}(\eta)$ indicate the maps
 $\theta\mapsto\ddot{\ell}(t,\theta,\eta)$ and
$\eta\mapsto\ell_{t,\theta}(t,\theta,\eta)$, respectively. For
brevity, we denote
$\dot{\ell}_{0}=\dot{\ell}(\theta_{0},\theta_{0},\eta_{0})$,
$\ddot{\ell}_{0}=\ddot{\ell}(\theta_{0},\theta_{0},\eta_{0})$ and
$\ell_{0}^{(3)}=\ell^{(3)}(\theta_{0},\theta_{0},\eta_{0})$, where
$\theta_0$, $\eta_0$ are the true values of $\theta$ and $\eta$.
Of course, we can write $\tilde{\ell}(X)$ as $\dot{\ell}_{0}(X)$.
$\|\cdot\|$ and $\|\cdot\|_{2}$ indicate the Euclidean norm and
$L_{2}$ norm, respectively. The notations $\apgt$ and $\aplt$ mean
greater than, or smaller than, up to a universal constant. The
symbols $\mathbb{P}_{n}$ and
$\mathbb{G}_n\equiv\sqrt{n}(\mathbb{P}_{n}-P)$ are used for the
empirical distribution and the empirical processes of the
observations, respectively.

\subsection{Main Assumptions}
We now make the following three classes of assumptions: Rate
assumptions (R1) for the penalized nuisance parameter and the
estimated penalty term; Smoothness assumptions (S1-S2) and
Empirical processes assumptions (E1) for $\ell(t,\theta,\eta)$ and
its related derivatives.

\begin{itemize}
\item[R1]: Assume:
\begin{eqnarray}
&&d(\hat{\eta}_{\tilde{\theta}_{n},\lambda_{n}},\eta_{0})=O_{P}(\lambda_{n}+\|\tilde{\theta}_{n}-\theta_{0}\|
)\label{penrate}\\
&&\mbox{and}\;\;\nonumber\\
&&\lambda_{n}J(\hat{\eta}_{\tilde{\theta}_{n},\lambda_{n}})=O_{P}(\lambda_{n}+\|\tilde{\theta}_{n}-\theta_{0}\|).\label{jrate}
\end{eqnarray}
\item[S1]: The maps
\begin{eqnarray}
(t,\theta,\eta)&\mapsto&\frac{\partial^{l+m}}{\partial
t^{l}\partial \theta^{m}}\ell(t,\theta,\eta)\label{mapcon}
\end{eqnarray}
have integrable envelope functions in $L_{1}(P)$ in some
neighborhood of $(\theta_{0},\theta_{0},\eta_{0})$, for
$(l,m)=(0,0), (1,0), (2,0), (3,0), (1,1), (1,2), (2,1)$.
\item[S2]: Assume:
\begin{eqnarray}
&&P\ddot{\ell}(\theta_{0},\theta_{0},\eta)
-P\ddot{\ell}(\theta_{0},\theta_{0},\eta_{0})=
O(d(\eta,\eta_{0}))\label{smcon2},\\
&&P\ell_{t,\theta}(\theta_{0},\theta_{0},\eta)
-P\ell_{t,\theta}(\theta_{0},\theta_{0},\eta_{0})
=O(d(\eta,\eta_{0}))\label{smcon3},\\
&&P\dot{\ell}(\theta_{0},\theta_{0},\eta)
=O(d^{2}(\eta,\eta_{0})),\label{nobias3}
\end{eqnarray}
for all $\eta$ in
some neighborhood of $\eta_{0}$.
\item[E1]: For all random
sequences $\tilde{\theta}_{n}=\hat{\theta}_{n}+o_{P}(1)$ and
$\bar{\theta}_{n}=\theta_{0}+o_{P}(1)$, we have
\begin{eqnarray}
\mathbb{G}_{n}(\dot{\ell}(\theta_{0},\theta_{0},\hat{\eta}_{\tilde{\theta}_{n},\lambda_{n}})-\dot{\ell}_{0})
&=&O_{P}(n^{\frac{1}{4k+2}}(\lambda_{n}+\|\tilde{\theta}_{n}-\theta_{0}\|))\label{smcon1},\\
\mathbb{G}_{n}(\ddot{\ell}(\theta_{0},\tilde{\theta}_{n},\hat{\eta}_{\tilde{\theta}_{n},\lambda_{n}}))&=&O_{P}(1),\label{eg2em1}\\
\mathbb{G}_{n}(\ell_{t,\theta}(\theta_{0},\bar{\theta}_{n},\hat{\eta}_{\tilde{\theta}_{n},\lambda_{n}}))&=&O_{P}(1),\label{eg2em2}\\
(\mathbb{P}_{n}-P)\ell^{(3)}(\bar{\theta}_{n},\tilde{\theta}_{n},\hat{\eta}_{\tilde{\theta}_{n},\lambda_{n}})&=&o_{P}(1).\label{eg2em3}
\end{eqnarray}
\end{itemize}

Assumption R1 implicitly assumes that we have a metric or
topology defined on the set of possible values of the nuisance
parameter $\eta$. The form of $d(\eta,\eta_{0})$ may vary for
different situations and does not need to be specified in this
subsection beyond the given conditions.
(\ref{penrate}) implies that
$\hat{\eta}_{\tilde{\theta}_{n},\lambda_{n}}$ is consistent for
$\eta_{0}$ as $\tilde{\theta}_{n}\rightarrow\theta_{0}$ in
probability. Additionally, from (\ref{jrate}) we know that the
smoothing parameter $\lambda_{n}$ plays a role in determining the
complexity degree of the estimated nuisance parameter.
(\ref{jrate}) implies that $J(\hat{\eta}_{\lambda_{n}})=O_{P}(1)$
if the $\hat{\theta}_{\lambda_{n}}$ is asymptotically normal,
which has been shown in (\ref{penmle}). Note that
$J(\hat{\eta}_{\tilde{\theta}_{n},0})\geq
J(\hat{\eta}_{\tilde{\theta}_{n},\lambda_{n}})$, where
$\hat{\eta}_{\theta,0}=\hat{\eta}_{\theta}\equiv
argmax_{\eta\in\mathcal{H}}\log lik(\theta,\eta)$ for a fixed
$\theta$, based on the inequality that $\log
lik_{\lambda_{n}}(\tilde{\theta}_{n},
\hat{\eta}_{\tilde{\theta}_{n},0})\leq\log
lik_{\lambda_{n}}(\tilde{\theta}_{n},\hat{\eta}_{\tilde{\theta}_{n},\lambda_{n}})$.

Clearly, the assumptions S1 and S2 are separately the smoothness
conditions for the Euclidean parameters $(t,\theta)$ and the
infinite dimensional nuisance parameter $\eta$. The
boundedness of the Fr\'{e}chet derivatives of the maps
$\eta\mapsto\ddot{\ell}(\theta_{0},\theta_{0},\eta)$ and
$\eta\mapsto\ell_{t,\theta}(\theta_{0},\theta_{0},\eta)$ ensures
the validity of conditions (\ref{smcon2}) and (\ref{smcon3}).
Based on the discussions in section 2 of \cite{ck062}, under the
given regularity conditions, it suffices to show (\ref{nobias3})
if the map $\eta\mapsto\dot{\ell}(\theta_{0},\theta_{0},\eta)$ is
Fr\'{e}chet differentiable and the map $\eta\mapsto
lik(\theta_{0},\eta)$ is second order Fr\'{e}chet differentiable.

Condition (\ref{smcon1}) is concerned with the asymptotic
equicontinuity of the empirical process measure of
$\dot{\ell}(\theta_{0},\theta_{0},\eta)$ with $\eta$ ranging
around the neighborhood of $\eta_{0}$. It suffices to show
(\ref{eg2em1}) and (\ref{eg2em2}) if
$\mathbb{G}_{n}(\ddot{\ell}(\theta_{0},\tilde{\theta}_{n},\hat{\eta}_{\tilde{\theta}_{n},\lambda_{n}})-\ddot{\ell}_{0})=o_{P}(1)$
and
$\mathbb{G}_{n}(\ell_{t,\theta}(\theta_{0},\bar{\theta}_{n},\hat{\eta}_{\tilde{\theta}_{n},\lambda_{n}})-\ell_{t,\theta}(\theta_{0},\theta_{0},\eta_{0}))
=o_{P}(1)$, provided $\ddot{\ell}_{0}$ and
$\ell_{t,\theta}(\theta_{0},\theta_{0},\eta_{0})$ are square
integrable. Thus we will be able to use technical tools T2 and
T6 given in the
appendix to show (\ref{smcon1})--(\ref{eg2em2}). For the
verification of (\ref{eg2em3}), we need to make use of a
Glivenko-Cantelli theorem for classes of functions that change
with $n$ which is a modification of theorem 2.4.3 in \cite{vw96}
and is explained in the appendix.

In principle, assumptions S1, S2 and E1 on the functions
of the least favorable submodel directly imply the following
empirical no-bias conditions:
\begin{eqnarray}
\pp_{n}\dot{\ell}(\theta_{0},\tilde{\theta}_{n},\hat{\eta}_{\tilde{\theta}_{n},\lambda_{n}})&=&\pp_{n}\tilde{\ell}_{0}+
O_{P}(\lambda_{n}+\|\tilde{\theta}_{n}
-\theta_{0}\|)^{2},\label{emno1}\\
\mathbb{P}_{n}\ddot{\ell}(\theta_{0},\tilde{\theta}_{n},
\hat{\eta}_{\tilde{\theta}_{n},\lambda_{n}})&=&P\ddot{\ell}_{0}+
O_{P}(\lambda_{n}+\|\tilde{\theta}_{n}
-\theta_{0}\|).\label{emno2}
\end{eqnarray}
The derivations of (\ref{emno1}) and (\ref{emno2}) are simply
based on the regular Taylor expansions around the true values. The
detailed arguments can be found in the proof of
lemmas~1 and~2 in~\cite{ck062}. The two empirical no-bias
conditions ensure that the penalized profile likelihood behaves
like a penalized likelihood in the parametric model
asymptotically and therefore yields a second order asymptotic
expansion of the penalized profile log-likelihood.

\section{The Penalized Convergence Rate}\label{rate}
In the previous section, we have imposed two assumptions about the
convergence rates of the estimated nuisance parameter and the
order of the estimated penalty term, i.e. (\ref{penrate}) and
(\ref{jrate}). To compute the convergence rates, we present three
different theorems below which require different sets of
conditions. These theorems can be viewed as extension of
general results on M-estimators to penalized M-estimators,
and are therefore of independent interest.
We first state the classical definitions for the
covering number (entropy number) and bracketing number (bracketing
entropy number) for a class of functions.

{\it Definition:} Let $\mathcal{A}$ be a subset of a (pseudo-)
metric space $(\mathcal{L},d)$ of real-valued functions. The
$\delta$-covering number $N(\delta,\mathcal{A},d)$ of
$\mathcal{A}$ is the smallest $N$ for which there exist functions
$a_{1},\ldots,a_{N}$ in $\mathcal{L}$, such that for each
$a\in\mathcal{A}$, $d(a,a_{j})\leq\delta$ for some
$j\in\{1,\ldots,N\}$. The $\delta$-bracketing number
$N_{B}(\delta,\mathcal{A},d)$ is the smallest $N$ for which there
exist pairs of functions
$\{[a_{j}^{L},a_{j}^{U}]\}_{j=1}^{N}\subset\mathcal{L}$, with
$d(a_{j}^{L},a_{j}^{U})\leq\delta$, $j=1,\ldots,N$, such that for
each $a\in\mathcal{A}$ there is a $j\in\{1,\ldots,N\}$ such that
$a_{j}^{L}\leq a\leq a_{j}^{U}$. The $\delta$-entropy number
($\delta$-bracketing entropy number) is defined as
$H(\delta,\mathcal{A},d)=\log N(\delta,\mathcal{A},d)$
($H_{B}(\delta,\mathcal{A},d)=\log N_{B}(\delta,\mathcal{A},d)$).

Before we present the first theorem, define
\begin{eqnarray*}
\mathcal{K}=\left\{\frac{\ell_{\theta,\eta}(X)-\ell_{0}(X)}{1+J(\eta)}:\|\theta-\theta_{0}\|\leq
C_{1}, \|\eta-\eta_{0}\|_{\infty}\leq C_{1},
J(\eta)<\infty\right\},
\end{eqnarray*}
for a known constant $C_1<\infty$:
\begin{theorem}\label{ratethm1}
Assume conditions (\ref{entrocon}), (\ref{eg1bou}),
(\ref{egsm3}) and (\ref{egsm1}) below hold for every
$\theta\in\Theta_{n}$ and $\eta\in\mathcal{V}_{n}$:
\begin{eqnarray}
&&H_{B}(\epsilon,\mathcal{K},L_{2}(P))\aplt\epsilon^{-1/k}\label{entrocon},\\
&&p_{\theta,\eta}/p_{\theta,\eta_{0}}\;\;\mbox{is bounded away from zero and infinity},\label{eg1bou}\\
&&\|\ell_{\theta,\eta}-\ell_{0}\|_{2}\aplt\|\theta-\theta_{0}\|+d_{\theta}(\eta,\eta_{0}),\label{egsm3}\\
&&P(\ell_{\theta,\eta}-\ell_{\theta,\eta_{0}})\aplt-d_{\theta}^{2}(\eta,\eta_{0})+\|\theta-\theta_{0}\|^{2}.\label{egsm1}
\end{eqnarray}
Then we have
\begin{eqnarray*}
d_{\tilde{\theta}_{n}}(\hat{\eta}_{\tilde{\theta}_{n},\lambda_{n}},\eta_{0})=O_{P}(\lambda_{n}+\|\tilde{\theta}_{n}-\theta_{0}\|),\\
\lambda_{n}J(\hat{\eta}_{\tilde{\theta}_{n},\lambda_{n}})=O_{P}(\lambda_{n}+\|\tilde{\theta}_{n}-\theta_{0}\|),
\end{eqnarray*}
for
$(\tilde{\theta}_{n},\hat{\eta}_{\tilde{\theta}_{n},\lambda_{n}})$
satisfying $P(\tilde{\theta}_{n}\in\Theta_{n},
\hat{\eta}_{\tilde{\theta}_{n},\lambda_{n}}\in\mathcal{V}_{n})\rightarrow
1$.
\end{theorem}

Condition (\ref{entrocon}) determines the order of the
increments of the empirical processes indexed by
$\ell_{\theta,\eta}$. A detailed discussion about how to
compute the increments of the empirical processes can be found in
chapter~5 of \cite{vg00}. Condition (\ref{eg1bou}) is
equivalent to the condition that $p_{\theta,\eta}$ is bounded away
from zero uniformly in $x$ for $(\theta,\eta)$ ranging over
$\Theta_{n}\times\mathcal{V}_{n}$. Given that the distance function
$d_{\theta}(\eta,\eta_{0})$ in (\ref{egsm3}) is just
$\|p_{\theta,\eta}-p_{0}\|_{2}$, (\ref{egsm3}) trivially holds
provided that condition (\ref{eg1bou}) holds. For the
verification of (\ref{egsm1}), we can do an analysis as follows. The
natural Taylor expansions of the criterion function
$(\theta,\eta)\mapsto P\ell_{\theta,\eta}$ around the maximum
point $(\theta_{0},\eta_{0})$ implies that
$P(\ell_{\theta,\eta_{0}}-\ell_{\theta_{0},\eta_{0}})\apgt
-\|\theta-\theta_{0}\|^{2}$, and (\ref{disine}) implies that
$P(\ell_{\theta,\eta}-\ell_{0})\leq-\int(\sqrt{p_{\theta,\eta}}-\sqrt{p_{0}})^{2}d\mu\leq-\|p_{\theta,\eta}-p_{0}\|_{2}^{2}$
given condition (\ref{eg1bou}).

We now apply theorem~\ref{ratethm1} to derive the related
convergence rates in the partly linear model in
corollary~\ref{eg1rate}. However, we need to strengthen our
previous assumptions to require the existence of a known
$M<\infty$ such that $\eta\in\mathcal{H}_{k}^{M}$, where
$\mathcal{H}_{k}^{M}=\mathcal{H}_{k}\cap\{\|\eta\|_{\infty}\leq
M\}$ and that the density for the joint distribution $(U,V,C)$ is
strictly positive and finite. The additional assumptions here
guarantee condition (\ref{eg1bou}). The following
theorem~\ref{nuirate} and theorem~\ref{ratethm3} can also be
employed to derive the convergence rate of the non-penalized
estimated nuisance parameter by setting $\lambda_{n}$ to zero.
However, we would need to assume that
$f\in\{g:\|g\|_{\infty}+J(g)\leq
\tilde{M}\}$ for some known $\tilde{M}$ when applying these
theorems. Thus we can argue that the the penalized method enables
a relaxation of the assumptions needed for the nuisance parameter.

\begin{corollary}\label{eg1rate}
Under the above set-up for the partly linear normal model with
current status data, we have, for $\tilde{\theta}_n=\theta_0+
o_P(1)$,
\begin{eqnarray}
\|\hat{f}_{\tilde{\theta}_{n},\lambda_{n}}-f_{0}\|_{2}&=&O_{P}(\lambda_{n}+\|\tilde{\theta}_{n}-\theta_{0}\|),\label{eg1pratre}\\
\lambda_{n}J(\hat{f}_{\tilde{\theta}_{n},\lambda_{n}})&=&O_{P}(\lambda_{n}+\|\tilde{\theta}_{n}-\theta_{0}\|).\label{eg1j}
\end{eqnarray}
Moreover,
if we also assume that $f\in\{g:\|g\|_{\infty}+J(g)\leq \tilde{M}\}$
for some known $\tilde{M}$, then
\begin{eqnarray}
\|\hat{f}_{\tilde{\theta}_{n}}-f_{0}\|_{2}&=&O_{P}(n^{-k/(2k+1)}+\|\tilde{\theta}_{n}-\theta_{0}\|),\label{eg10ratre}
\end{eqnarray}
provided condition (\ref{smooth}) holds.
\end{corollary}

\begin{remark}
Corollary~\ref{eg1rate} implies that the convergence
rate of the estimated nuisance parameter is slower than that of
the the regular nuisance parameter by comparing (\ref{eg1pratre})
and (\ref{eg10ratre}). This result is not surprising since the
slower rate is the trade off for the smoother nuisance parameter
estimator. However the advantage of the penalized profile sampler
is that we can control the convergence rate by assigning the
smoothing parameter with different rates.
Corollary~\ref{eg1rate} also indicates that
$\|\hat{f}_{\lambda_{n}}-f_{0}\|_{2}=O_{P}(\lambda_{n})$ and
$\|\hat{f}_{n}-f_{0}\|_{2}=O_{P}(n^{-k/(k+2)})$. Note that the
convergence rate of the maximum penalized likelihood estimator,
$O_{P}(\lambda_{n})$, is deemed as the optimal rate in \cite{w90}.
Similar remarks also hold for corollary~\ref{eg2rate} below.
\end{remark}

The boundedness condition (\ref{eg1bou}) appears hard to achieve
in some examples. Hence we propose theorem~\ref{nuirate} below to
relax this condition by choosing the criterion function
$m_{\theta,\eta}=\log[(p_{\theta,\eta}+p_{\theta,\eta_{0}})/2p_{\theta,\eta_{0}}]$.
Obviously, $m_{\theta,\eta}$ is trivially bounded away from zero.
It is also bounded above for $(\theta,\eta)$ around the their
true values if $p_{\theta,\eta_{0}}(x)$ is bounded away from zero
uniformly in $x$ and $p_{\theta,\eta}$ is bounded above. The first
condition is satisfied if the map $\theta\mapsto
p_{\theta,\eta_{0}}(x)$ is continuous around $\theta_{0}$ and
$p_{0}(x)$ is uniformly bounded away from zero. The second
condition is trivially satisfied in the semiparametric logistic
regression model by the given form of the density. The boundedness of
$m_{\theta,\eta}$ thus permits the application of
lemma~\ref{lerate} below which is used to verify condition
(\ref{ratea2}) in the following theorem:
\begin{theorem}\label{nuirate}
Assume for any given $\theta\in\Theta_{n}$, $\hat{\eta}_{\theta}$
satisfies
$\mathbb{P}_{n}m_{\theta,\hat{\eta}_{\theta}}\geq\mathbb{P}_{n}m_{\theta,\eta_{0}}$
for given measurable functions $x\mapsto m_{\theta,\eta}(x)$.
Assume conditions (\ref{ratea1}) and (\ref{ratea2}) below hold for
every $\theta\in\Theta_{n}$, every $\eta\in\mathcal{V}_{n}$ and
every $\epsilon>0$:
\begin{eqnarray}
&&P(m_{\theta,\eta}-m_{\theta,\eta_{0}})\aplt
-d^{2}_{\theta}(\eta,\eta_{0})+\|\theta-\theta_{0}\|^{2},\label{ratea1}\\
&&E^{\ast}\sup_{\theta\in\Theta_{n},\eta\in\mathcal{V}_{n},\|\theta-\theta_{0}\|<\epsilon,d_{\theta}(\eta,\eta_{0})<\epsilon}|\mathbb{G}_{n}
(m_{\theta,\eta}-m_{\theta,\eta_{0}})|\aplt\phi_{n}(\epsilon).\label{ratea2}
\end{eqnarray}
Suppose that (\ref{ratea2}) is valid for functions $\phi_{n}$ such
that $\delta\mapsto\phi_{n}(\delta)/\delta^{\alpha}$ is decreasing
for some $\alpha<2$ and sets $\Theta_{n}\times\mathcal{V}_{n}$
such that
$P(\tilde{\theta}\in\Theta_{n},\hat{\eta}_{\tilde{\theta}}\in\mathcal{V}_{n})\rightarrow
1$. Then
$d_{\tilde{\theta}}(\hat{\eta}_{\tilde{\theta}},\eta_{0})\leq
O_{P}^{\ast}(\delta_{n}+\|\tilde{\theta}-\theta_{0}\|)$ for any
sequence of positive numbers $\delta_{n}$ such that
$\phi_{n}(\delta_{n})\leq\sqrt{n}\delta_{n}^{2}$ for every $n$.
\end{theorem}

Lemma~\ref{lerate} below is presented to verify the modulus
condition for the continuity of the empirical process in
(\ref{ratea2}). Let $\mathcal{S}_{\delta}=\{x\mapsto
m_{\theta,\eta}(x)-m_{\theta,\eta_{0}}(x):
d_{\theta}(\eta,\eta_{0})<\delta, \|\theta-\theta_{0}\|<\delta\}$
and write
\begin{eqnarray}
K(\delta,\mathcal{S}_{\delta},L_{2}(P))=\int_{0}^{\delta}\sqrt{1+H_{B}(\epsilon,\mathcal{S}_{\delta},L_{2}(P))}d\epsilon:\label{thm2k}
\end{eqnarray}

\begin{lemma}\label{lerate}
Suppose the functions $(x,\theta,\eta)\mapsto
m_{\theta,\eta}(x)$ are uniformly bounded for $(\theta,\eta)$
ranging over a neighborhood of $(\theta_{0},\eta_{0})$ and that
\begin{eqnarray*}
P(m_{\theta,\eta}-m_{\theta_{0},\eta_{0}})^{2}\aplt
d_{\theta}^{2}(\eta,\eta_{0})+\|\theta-\theta_{0}\|^{2}.
\end{eqnarray*}
Then condition (\ref{ratea2}) is satisfied for any functions
$\phi_{n}$ such that
\begin{eqnarray*}
\phi_{n}(\delta)\geq
K(\delta,\mathcal{S}_{\delta},L_{2}(P))\left(1+\frac{K(\delta,\mathcal{S}_{\delta},L_{2}(P))}{\delta^{2}\sqrt{n}}\right)
\end{eqnarray*}
Consequently, in the conclusion of the above theorem we may use
$K(\delta,\mathcal{S}_{\delta},L_{2}(P))$ rather than
$\phi_{n}(\delta)$.
\end{lemma}

\begin{remark}
Theorem~\ref{eg2rate} and lemma~\ref{lerate} are theorem
3.2 and lemma 3.3 in \cite{mv99}, respectively. We can apply
theorem~\ref{nuirate} to the penalized semiparametric logistic
regression model by including $\lambda$ in $\theta$, i.e.
$m_{\theta,\lambda,\eta}=m_{\theta,\eta}-\frac{1}{2}\lambda^{2}(J^{2}(\eta)-
J^{2}(\eta_{0}))$. This is accomplished in the following
corollary. Note that we assume that the uniform norm and
Sobolev norm of $\eta$ are bounded above with known upper bounds
when deriving (\ref{eg1ratre}) of the corollary, but this
assumption is not needed for~(\ref{eg2pratre}) and~(\ref{eg2j}).
\end{remark}

\begin{corollary}\label{eg2rate}
Under the above set-up for the semiparametric logistic regression
model, we have for $\lambda_n$ satisfying condition~(\ref{smooth})
and any $\tilde{\theta}_{n}\overset{p}{\rightarrow}\theta_{0}$ that
\begin{eqnarray}
\|\hat{\eta}_{\tilde{\theta}_{n},\lambda_{n}}-\eta_{0}\|_{2}&=&O_{P}(\lambda_{n}+\|\tilde{\theta}_{n}-\theta_{0}\|),\label{eg2pratre}\\
\lambda_{n}J(\hat{\eta}_{\tilde{\theta}_{n},\lambda_{n}})&=&O_{P}(\lambda_{n}+\|\tilde{\theta}_{n}-\theta_{0}\|).\label{eg2j}
\end{eqnarray}
If we also
assume that $\eta\in\{g:\|g\|_{\infty}+J(g)\leq \tilde{M}\}$
for some known $\tilde{M}$, then
\begin{eqnarray}
\|\hat{\eta}_{\tilde{\theta}_{n}}-\eta_{0}\|_{2}&=&O_{P}(n^{-k/(2k+1)}+\|\tilde{\theta}_{n}-\theta_{0}\|).\label{eg1ratre}
\end{eqnarray}
\end{corollary}

\begin{remark}
Corollary~\ref{eg1rate} and~\ref{eg2rate} imply that
$J(\hat{\eta}_{\lambda_{n}})=O_{P}(1)$ and
$J(\hat{f}_{\lambda_{n}})=O_{P}(1)$, respectively. Thus the
maximum likelihood estimators of the nuisance parameters in the
two examples of this paper are consistent in the uniform norm,
i.e. $\|\hat{\eta}_{\lambda_{n}}-\eta_{0}\|_{\infty}=o_{P}(1)$ and
$\|\hat{f}_{\lambda_{n}}-f_{0}\|_{\infty}=o_{P}(1)$, since the
sequences $\hat{\eta}_{n}$ and $\hat{f}_{n}$ consist of smooth
functions defined on a compact set with asymptotically bounded
first-order derivatives.
\end{remark}

The preceding two theorems imply that the convergence rate of the
penalized estimated nuisance parameter is affected by the assigned
smoothness parameter. However, the next theorem shows that, under
different conditions, the above phenomena may not hold. Let
\begin{eqnarray*}
l^{\lambda_{n}}_{\theta,\eta,h}&=&\frac{\partial}{\partial
t}|_{t=0}\log
lik_{\lambda_{n}}(\theta,\eta_{t})=A_{\theta,\eta}h-2\lambda_{n}^{2}\int
h^{(k)}\eta^{(k)}dz,\\
V(\theta,\eta)h&=&PA_{\theta,\eta}h,\\
\mathbb{V}_{n}(\theta,\eta)h&=&\mathbb{P}_{n}A_{\theta,\eta}h,
\end{eqnarray*}
where $\eta_{t}=\eta+th$ for $h\in\mathcal{H}_{k}$ and
$A_{\theta,\eta}$ is the appropriate score operator for the
model. Note that
$\eta_{t}\in\mathcal{H}_{k}$ for sufficiently small $t$. Obviously
$\mathbb{P}_{n}l^{\lambda_{n}}_{\tilde{\theta}_{n},\hat{\eta}_{\tilde{\theta}_{n},\lambda_{n}},h}=0$
and $V(\theta_{0},\eta_{0})h=0$. We assume that the maps
$h\mapsto V(\theta,\eta)h$ and $h\mapsto
\mathbb{V}_{n}(\theta,\eta)h$ are uniformly bounded such that
$\mathbb{V}_{n}$ and $V$ can be viewed as maps from the parameters
set $\Theta\times\mathcal{H}_{k}$ into
$\ell^{\infty}(\mathcal{H}_{k})$. Further we require the following
regularity conditions: For some $C_{2}>0$,
\begin{eqnarray}
&&\left\{A_{\theta,\eta}h:\|\theta-\theta_{0}\|<C_{2},d_{\theta}(\eta,\eta_{0})<C_{2},h\in\mathcal{H}_{k}\right\}\;\;\mbox{is
$P$-Donsker},\label{con1th3}\\
&&\sup_{h\in\mathcal{H}_{k}}P(A_{\theta,\eta}h-A_{\theta_{0},\eta_{0}}h)^{2}\rightarrow
0,\;\mbox{as}\;
\theta\rightarrow\theta_{0}\;\;\mbox{and}\;\;\eta\rightarrow\eta_{0}.\label{con2th3}
\end{eqnarray}

\begin{theorem}\label{ratethm3}
Suppose that
$V(\cdot,\cdot):\Theta\times\mathcal{H}_{k}\mapsto\ell^{\infty}(\mathcal{H}_{k})$
is Fr\'{e}chet differentiable at $(\theta_{0},\eta_{0})$ with
derivative $\dot{V}(\cdot,\cdot):\mathbb{R}^{d}\times
lin\mathcal{H}_{k}\mapsto\ell^{\infty}(\mathcal{H}_{k})$ such that
the map
$\dot{V}(0,\cdot):lin\mathcal{H}_{k}\mapsto\ell^{\infty}(\mathcal{H}_{k})$
is invertible with an inverse that is continuous on its range.
Furthermore, we assume that (\ref{con1th3}) and (\ref{con2th3})
hold. Then
\begin{eqnarray}
d_{\tilde{\theta}_{n}}(\hat{\eta}_{\tilde
{\theta}_{n},\lambda_{n}},\eta_{0})=O_{P}(n^{-1/2}+\|\tilde{\theta}_{n}-\theta_{0}\|+\lambda_{n}^{2}J^{2}(\hat{\eta}_{\tilde{\theta}_{n}
,\lambda_{n}})),\label{thm3rate}
\end{eqnarray}
for $\tilde{\theta}_{n}\rightarrow \theta_{0}$ and
$\hat{\eta}_{\tilde{\theta}_{n},\lambda_{n}}\rightarrow \eta_{0}$
in probability.
\end{theorem}

\begin{remark}
The preceding theorem is a variation of theorems used in
\cite{m95} and \cite{v94}, among others, to prove the asymptotic
normality of the maximum likelihood estimator
$(\hat{\theta}_{n},\hat{\eta}_{n})$. If we can show that
$\lambda_{n}J(\hat{\eta}_{\tilde{\theta}_{n}
,\lambda_{n}})=O_{P}(\lambda_{n}+\|\tilde{\theta}_{n}-\theta_{0}\|)$
by some other means, then (\ref{thm3rate}) implies that
$d_{\tilde{\theta}_{n}}(\hat{\eta}_{\tilde
{\theta}_{n},\lambda_{n}},\eta_{0})=O_{P}(n^{-1/2}+\|\tilde{\theta}_{n}-\theta_{0}\|)$.
This indicates that the smoothing effect of the penalized method
does not occur, which may be due to some very smooth non-penalized
estimated nuisance parameter. The high degree of the smoothness of
the non-penalized estimated nuisance parameter can be deduced from
its fast convergence rate which equals the parametric rate in
this instance.
\end{remark}

\section{Main Results and Implications}
In this section we first present second order asymptotic expansion
of the log-profile penalized likelihood which prepare us for
deriving the main results about the higher order structure of the
penalized profile sampler. The assumptions in section 3 and
condition~(\ref{smooth}) are assumed throughout.

\begin{theorem}\label{penexp}
Given $\tilde{\theta}_{n}=\hat{\theta}_{\lambda_{n}}+o_{P}(1)$,
we have
\begin{eqnarray}
\sqrt{n}(\hat{\theta}_{\lambda_{n}}-\theta_{0})&=&\frac{1}{\sqrt{n}}\sum_{i=1}^{n}\tilde{I}_{0}^{-1}\tilde{\ell}_{0}(X_{i})+
O_{P}(n^{1/2}\lambda_{n}^{2}),\label{penmle}\\
\log pl_{\lambda_{n}}(\tilde{\theta}_{n})&=&\log
pl_{\lambda_{n}}(\hat{\theta}_{\lambda_{n}})-\frac{n}{2}(\tilde{\theta}_{n}-
\hat{\theta}_{\lambda_{n}})^{T}\tilde{I}_{0}(\tilde{\theta}_{n}-
\hat{\theta}_{\lambda_{n}})\label{penexpf}\\&+&O_{P}(g_{\lambda_{n}}(\|\tilde{\theta}_{n}-\hat{\theta}_{\lambda_{n}}\|))
\nonumber,
\end{eqnarray}
where
$g_{\lambda_{n}}(w)=nw^{3}+nw^{2}\lambda_{n}+nw\lambda_{n}^{2}+n^{1/2}\lambda_{n}^{2}$,
provided the efficient information $\tilde{I}_{0}$ is positive
definite.
\end{theorem}

\begin{remark}
The results in theorem~\ref{penexp} are useful in there own right
for inference about $\theta$. (\ref{penmle}) is a second higher
order frequentist result in penalized semiparametric
estimation regarding the asymptotic linearity of
the maximum penalized likelihood estimator of $\theta$.
\end{remark}

We now state the main results on the penalized posterior profile
distribution. A preliminary result, theorem~\ref{pdist} with
corollary~\ref{pmome} below, shows that the penalized posterior
profile distribution is asymptotically close enough to the
distribution of a normal random variable with mean
$\hat{\theta}_{\lambda_{n}}$ and variance $(n\tilde{I}_{0})^{-1}$
with second order accuracy, which is controlled by the smoothing
parameter. Similar conclusions also hold for the penalized
posterior moments. Another main result, theorem~\ref{pquant},
shows that the penalized posterior profile log-likelihood can be
used to achieve second order accurate frequentist inference for
$\theta$.

Let $\tilde{P}_{\theta|\tilde{X}}^{\lambda_{n}}$ be the penalized
posterior profile distribution of $\theta$ with respect to the prior
$\rho(\theta)$. Define
\begin{eqnarray*}
\Delta_{\lambda_{n}}(\theta)&\equiv&n^{-1}\{\log
pl_{\lambda_{n}}(\theta)-\log
pl_{\lambda_{n}}(\hat{\theta}_{\lambda_{n}
})\}\\
&=&n^{-1}(\ell_{n}(\theta,\hat{\eta}_{\theta,\lambda_{n}})-\ell_{n}(\hat{\theta}_{\lambda_{n}},\hat{\eta}_{\lambda_{n}}))-
n^{-1}\lambda_{n}^{2}(J^{2}(\hat{\eta}_{\theta,\lambda_{n}})-
J^{2}(\hat{\eta}_{\lambda_{n}})).
\end{eqnarray*}

\begin{theorem}\label{pdist}
Assume that
\begin{eqnarray}
\Delta_{\lambda_{n}}(\tilde{\theta}_{n})=o_P(1) \;\;\mbox{ implies
}\;\; \tilde{\theta}_{n}=\theta_{0}+o_P(1),\label{asyuni}
\end{eqnarray}
for every random $\left\{\tilde{\theta}_{n}\right\}\in\Theta$. If
$\rho(\theta_{0})>0$ and $\rho(\cdot)$ has continuous and finite
first order derivative in some neighborhood of $\theta_{0}$, then
we have, for any $-\infty<\xi<\infty$,
\begin{eqnarray}
\sup_{\xi\in\mathbb{R}^{d}}\left|\tilde{P}_{\theta|\tilde{X}}^{\lambda_{n}}(\sqrt{n}\tilde{I}_{0}^{1/2}(\theta-\hat{\theta}_{\lambda_{n}})
\leq\xi)-\Phi_{d}(\xi)\right|=
O_P(n^{1/2}\lambda_{n}^{2}),\label{ppostexp2}
\end{eqnarray}
where $\Phi_{d}(\cdot)$ is the distribution of the $d$-dimensional
standard normal random variable.
\end{theorem}

\begin{corollary}\label{pmome}
Under the assumptions of theorem~\ref{pdist}, we have that if
$\theta$ has finite second absolute moment, then
\begin{eqnarray}
\hat{\theta}_{\lambda_{n}}&=&E^{\lambda_{n}}_{\theta|\tilde{X}}(\theta)+O_{P}(\lambda_{n}^{2}),\label{pmon1}\\
\tilde{I}_{0}&=&n^{-1}(Var^{\lambda_{n}}_{\theta|\tilde{X}}(\theta))^{-1}+O_{P}(n^{1/2}\lambda_{n}^{2}),\label{pmon2}
\end{eqnarray}
where $E^{\lambda_{n}}_{\theta|\tilde{X}}(\theta)$ and
$Var^{\lambda_{n}}_{\theta|\tilde{X}}(\theta)$ are the penalized
posterior profile mean and penalized posterior profile covariance
matrix, respectively.
\end{corollary}

We now present another second order asymptotic frequentist
property of the penalized profile sampler in terms of quantiles.
The $\alpha$-th quantile of the penalized posterior profile
distribution, $\tau_{n\alpha}$, is defined as
$\tau_{n\alpha}=inf\lbrace \xi:
\tilde{P}^{\lambda_{n}}_{\theta|\tilde{X}}(\theta\leq
\xi)\geq\alpha\rbrace$. Without loss of generality,
$\tilde{P}^{\lambda_{n}}_{\theta|\tilde{X}}(\theta\leq
\tau_{n\alpha})=\alpha$. We can also define
$\kappa_{n\alpha}\equiv\sqrt{n}(\tau_{n\alpha}-\hat{\theta}_{\lambda_{n}})$,
i.e.,
$\tilde{P}^{\lambda_{n}}_{\theta|\tilde{X}}(\sqrt{n}(\theta-\hat{\theta}_{\lambda_{n}})\leq\kappa_{n\alpha})=\alpha$.

\begin{theorem}\label{pquant}
Under the assumptions of theorem~\ref{pdist} and assuming that
$\tilde{\ell}_{0}(X)$ has finite third moment with a nondegenerate
distribution, then there exists a $\hat{\kappa}_{n\alpha}$
based on the data such that
$P(\sqrt{n}(\hat{\theta}_{\lambda_{n}}-\theta_{0})\leq\hat{\kappa}_{n\alpha})=\alpha$
and
$\hat{\kappa}_{n\alpha}-\kappa_{n\alpha}=O_P(n^{1/2}\lambda_{n}^{2})$
for each choice of $\kappa_{n\alpha}$.
\end{theorem}

\begin{remark}
Theorem~\ref{pquant} ensures that there exists a unique
$\alpha$-th quantile for $\theta$ up to $O_{P}(\lambda_{n}^{2})$
in the frequentist set-up for each fixed $\tau_{n\alpha}$. Note
that $\tau_{n\alpha}$ is not unique if the dimension of $\theta$
is larger than one.
\end{remark}
\begin{remark}
Theorem~\ref{pdist}, corollary~\ref{pmome} and
theorem~\ref{pquant} above show that the penalized profile
sampler generates second order asymptotic frequentist valid
results in terms of distributions, moments and quantiles.
Moreover, the second order accuracy of this procedure is
controlled by the smoothing parameter.
\end{remark}

\begin{remark}
Another interpretation for the role of $\lambda_{n}$ in the
penalized profile sampler is that we can view $\lambda_{n}$ as the
prior on $J(\eta)$, or on $\eta$ to some extent. To see this, we
can write $lik_{\lambda_{n}}(\theta,\eta)$ in the following form:
\begin{eqnarray*}
lik_{\lambda_{n}}(\theta,\eta)=lik_{n}(\theta,\eta)\times\exp\left[-\frac{J^{2}(\eta)}{2(\frac{1}{2\lambda_{n}^{2}})}\right]
\end{eqnarray*}
This idea can be traced back to \cite{w90}. In other words, the
prior on $J(\eta)$ is a normal distribution with mean zero and
variance $(2\lambda_{n}^{2})^{-1}$. Hence it is natural to expect
$\lambda_{n}$ has some effect on the convergence rate of $\eta$.
Other possible priors on the functional parameter include
Dirichlet and Gaussian processes which are more
commonly used in nonparametric Bayesian methodology.
\end{remark}

\section{Examples (Continued)}
We now illustrate verification of the assumptions in
section~3.3 with the two example that were introduced in section~2.
Thus this section is a continuation of the earlier examples.

\subsection{Partly Linear Normal Model with Current Status Data}
We will concentrate on the estimation of the regression
coefficient $\theta$, considering the infinite dimensional
parameter $f\in\mathcal{H}_{k}^{M}$ as a nuisance parameter. The
score function of $\theta$, $\dot{\ell}_{\theta,f}$, is given as
follows:
\begin{eqnarray*}
\dot{\ell}_{\theta,f}(x)=uQ(x;\theta,f),
\end{eqnarray*}
where
\begin{eqnarray*}
Q(X;\theta,f)=(1-\Delta)\frac{\phi(q_{\theta,f}(X))}
{1-\Phi(q_{\theta,f}(X))}-\Delta\frac{\phi(q_{\theta,f}(X))}{\Phi(q_{\theta,f}(X))},
\end{eqnarray*}
$q_{\theta,f}(x)=c-\theta u-f(v)$, and $\phi$ is the
density of a standard normal random variable. The
least favorable direction at the true parameter value is:
\begin{eqnarray*}
h_{0}(v)=\frac{E_{0}(UQ^{2}(X;\theta,f)|V=v)}{E_{0}(Q^{2}(X;\theta,f)|V=v)},
\end{eqnarray*}
where $E_{0}$ is the expectation relative to the true parameters.
The derivation of $\dot{\ell}_{\theta,f}$ and $h_{0}(\cdot)$ is given
in \cite{ck062}. Thus, the least favorable submodel
can be constructed as follows:
\begin{eqnarray}
\ell(t,\theta,f)=\log lik(t,f_{t}(\theta,f)),\label{eg2lfs}
\end{eqnarray}
where $f_{t}(\theta,f)=f+(\theta-t)h_{0}$. By differentiating
(\ref{eg2lfs}) with respect to $t$ or $\theta$, we can obtain the
maps assessed in assumption S1,
$(t,\theta,f)\mapsto(\partial^{l+m}/\partial t^{l}\partial
\theta^{m})\ell(t,\theta,f)$. The concrete forms of these maps are
given in \cite{ck062} which considers a more rigid model with a
known upper bound on the $L_2$ norm of the $k$th derivative. The
rate assumptions (\ref{penrate}) and (\ref{jrate}) have been
verified previously in corollary~\ref{eg1rate}. The remaining
assumptions are verified in the following two lemmas:

\begin{lemma}\label{eg2le0}
Under the above set-up for the partly linear normal model with current
status data, assumptions S1, S2 and E1 are satisfied.
\end{lemma}

\begin{lemma}\label{eg2le}
Under the above set-up for the partly linear normal model with current
status data, condition (\ref{asyuni}) is satisfied.
\end{lemma}

\subsection{Semiparametric Logistic Regression}
In the semiparametric logistic regression model, we can obtain the
score function for $\theta$ and $\eta$ by similar analysis
performed in the first example, i.e.
$\dot{\ell}_{\theta,\eta}(x)=(y-F(\theta w+\eta(z)))w$ and
$A_{\theta,\eta}h_{\theta,\eta}(x)=(y-F(\theta
w+\eta(z)))h_{\theta,\eta}(z)$ for $J(h)<\infty$. And the least
favorable direction at the true parameter is given in \cite{mv99}:
\begin{eqnarray*}
h_{0}(z)=\frac{P_{0}[W\dot{F}(\theta_{0}
W+\eta_{0}(Z))|Z=z]}{P_{0}[\dot{F}(\theta_{0}
W+\eta_{0}(Z))|Z=z]},
\end{eqnarray*}
where $\dot{F}(u)=F(u)(1-F(u))$. The above assumptions plus the
requirement that $J(h_{0})<\infty$ ensures the identifiability of
the parameters. Thus the least favorable submodel can be written
as:
\begin{eqnarray*}
\ell(t,\theta,\eta)=\log lik(t,\eta_{t}(\theta,\eta)),
\end{eqnarray*}
where $\eta_{t}(\theta,\eta)=\eta+(\theta-t)h_{0}$. By
differentiating $\ell(t,\theta,\eta)$ with respect to $t$ or
$\theta$, we obtain,
\begin{eqnarray*}
\dot{\ell}(t,\theta,\eta)&=&(y-F(tw+\eta(z)+(\theta-t)h_{0}(z)))(w-h_{0}(z)),\\
\ddot{\ell}(t,\theta,\eta)&=&-\dot{F}(tw+\eta(z)+(\theta-t)h_{0}(z))(w-h_{0}(z))^{2},\\
\ell_{t,\theta}(t,\theta,\eta)&=&-\dot{F}(tw+\eta(z)+(\theta-t)h_{0}(z))(w-h_{0}(z))h_{0}(z),\\
\ell^{(3)}(t,\theta,\eta)&=&-\ddot{F}(tw+\eta(z)+(\theta-t)h_{0}(z))(w-h_{0}(z))^{3},\\
\ell_{t,t,\theta}(t,\theta,\eta)&=&-\ddot{F}(tw+\eta(z)+(\theta-t)h_{0}(z))(w-h_{0}(z))^{2}h_{0}(z),\\
\ell_{t,\theta,\theta}(t,\theta,\eta)&=&-\ddot{F}(tw+\eta(z)+(\theta-t)h_{0}(z))(w-h_{0}(z))h_{0}^{2}(z),
\end{eqnarray*}
where $\ddot{F}(\cdot)$ is the second derivative of the function
$F(\cdot)$. The rate assumptions have been shown in
corollary~\ref{eg2rate}. The remaining assumptions are verified in
the following two lemmas:

\begin{lemma}\label{eg1le0}
Under the above set-up for the semiparametric logistic regression
model, assumptions S1, S2 and E1 are satisfied.
\end{lemma}

\begin{lemma}\label{eg1le}
Under the above set-up for the semiparametric logistic regression
model, condition (\ref{asyuni}) is satisfied.
\end{lemma}

\section{Future Work}
Our paper evaluates the penalized profile sampler method from the
frequentist view and discusses the effect of the smoothing parameter
on estimation accuracy. One potential problem of interest is
how to select a proper smoothing parameter in applications. A
formal study about the higher order comparisons between the
profile sampler procedure and fully Bayesian procedure \cite{s02},
which assign priors to both the finite dimensional parameter and the
infinite dimensional nuisance parameter, is also interesting. We
expect that the involvement of a suitable
prior on the infinite dimensional
parameter would at least not decrease the estimation accuracy of
the parameter of interest.

Another worthwhile avenue of research is
to develop analogs of the profile sampler and penalized profile
sampler to likelihood estimation under model misspecification and
to general M-estimation. Some first order results for this setting
in the case where the nuisance parameter may not be root-$n$
consistent have been developed for a weighted bootstrap
procedure in \cite{mk05}.

\section{Appendix}
We first present some technical tools about the entropy
calculations and increments of empirical processes which will be
employed in the proofs that follow.

{\indent \sc T1.} For each $0<C<\infty$ and $\delta>0$ we have
\begin{eqnarray}
H_{B}(\delta,\{\eta:\|\eta\|_{\infty}\leq C, J(\eta)\leq
C\},\|\cdot\|_{\infty})\aplt (\frac{C}{\delta})^{1/k},\label{sob1}\\
H(\delta,\{\eta:\|\eta\|_{\infty}\leq C, J(\eta)\leq
C\},\|\cdot\|_{\infty})\aplt (\frac{C}{\delta})^{1/k}.\label{sob2}
\end{eqnarray}

{\indent \sc T2.} Let $\mathcal{F}$ be a class of measurable
functions such that $Pf^{2}<\delta^{2}$ and $\|f\|_{\infty}\leq M$
for every $f$ in $\mathcal{F}$. Then
\begin{eqnarray*}
E_{P}^{\ast}\|\mathbb{G}_{n}\|_{\mathcal{F}}\aplt
K(\delta,\mathcal{F},L_{2}(P))\left(1+\frac{K(\delta,\mathcal{F},L_{2}(P))}{\delta^{2}\sqrt{n}}M\right),
\end{eqnarray*}
where
$K(\delta,\mathcal{F},\|\cdot\|)=\int_{0}^{\delta}\sqrt{1+H_{B}(\epsilon,\mathcal{F},\|\cdot\|)}d\epsilon$.

{\indent \sc T3.} Let $\mathcal{F}=\{f_{t}:t\in T\}$ be a class of
functions satisfying $|f_{s}(x)-f_{t}(x)|\leq d(s,t)F(x)$ for
every $s$ and $t$ and some fixed function $F$. Then, for any norm
$\|\cdot\|$,
\begin{eqnarray*}
N_{[]}(2\epsilon\|F\|,\mathcal{F},\|\cdot\|)\leq N(\epsilon,T,d).
\end{eqnarray*}

{\indent \sc T4.}
\begin{eqnarray}
-P_{\theta_{0}}\log\frac{p_{\theta}}{p_{\theta_{0}}}\geq\int\left(\sqrt{p_{\theta}}-\sqrt{p_{\theta_{0}}}\right)^{2}d\mu.\label{disine}
\end{eqnarray}

{\indent \sc T5.} Let $\mathcal{F}$ be a class of measurable
functions $f:\mathbb{D}\times\mathbb{W}\mapsto\mathbb{R}$ on a
product of a finite set and an arbitrary measurable space
$(\mathbb{W},\mathcal{W})$. Let $P$ be a probability measure on
$\mathbb{D}\times\mathbb{W}$ and let $P_{\mathcal{W}}$ be its
marginal on $\mathbb{W}$. For every $d\in\mathbb{D}$, let
$\mathcal{F}_{d}$ be the set of functions $w\mapsto f(d,w)$ as $f$
ranges over $\mathcal{F}$. If every class $\mathcal{F}_{d}$ is
$P$-Donsker with $\sup_{f\in\mathcal{F}}|Pf(d,W)|<\infty$ for
every $d$, then $\mathcal{F}$ is $P$-Donsker.

{\indent \sc T6.} Let $\mathcal{F}$ be a uniformly bounded class
of measurable functions such that for some measurable $f_{0}$,
$\sup_{f\in\mathcal{F}}\|f-f_{0}\|_{\infty}<\infty$. Moreover,
assume that $H_{B}(\epsilon,\mathcal{F},L_{2}(P))\leq
K\epsilon^{-\alpha}$ for some $K<\infty$ and $\alpha\in(0,2)$ and
for all $\epsilon>0$. Then
\begin{eqnarray*}
\sup_{f\in\mathcal{F}}\left[\frac{|(\mathbb{P}_{n}-P)(f-f_{0})|}{\|f-f_{0}\|_{2}^{1-\alpha/2}\vee
n^{(\alpha-2)/[2(2+\alpha)]}}\right]=O_{P}(n^{-1/2}).
\end{eqnarray*}

{\indent \sc T7.} For a probability measure $P$, let
$\mathcal{F}_{1}$ be a class of measurable functions
$f_{1}:\mathcal{X}\mapsto\mathbb{R}$, and let $\mathcal{F}_{2}$
denote a class of nondecreasing functions
$f_{2}:\mathbb{R}\mapsto[0,1]$ that are measurable for every
probability measure. Then,
\begin{eqnarray*}
H_{B}(\epsilon,\mathcal{F}_{2}(\mathcal{F}_{1}),L_{2}(P))\leq
2H_{B}(\epsilon/3,\mathcal{F}_{1},L_{2}(P))+\sup_{Q}H_{B}(\epsilon/3,\mathcal{F}_{2},L_{2}(Q)).
\end{eqnarray*}

{\indent \sc T8.} Let $\mathcal{F}$ and $\mathcal{G}$ be classes
of measurable functions. Then for any probability measure $Q$ and
any $1\leq r\leq\infty$,
\begin{eqnarray}
H_{B}(2\epsilon,\mathcal{F}+\mathcal{G},L_{r}(Q))\leq
H_{B}(\epsilon,\mathcal{F},L_{r}(Q))+
H_{B}(\epsilon,\mathcal{G},L_{r}(Q)),\label{product}
\end{eqnarray}
and, provided $\mathcal{F}$ and $\mathcal{G}$ are bounded by 1,
\begin{eqnarray}
H_{B}(2\epsilon,\mathcal{F}\times\mathcal{G},L_{r}(Q))\leq
H_{B}(\epsilon,\mathcal{F},L_{r}(Q))+
H_{B}(\epsilon,\mathcal{G},L_{r}(Q)).\label{times}
\end{eqnarray}

\begin{remark}\label{techdis}
The proof of T1 is found in \cite{bs67}. T1 implies that the
Sobolev class of functions with known bounded Sobolev norm is
$P$-Donsker. T2 and T3 are separately lemma 3.4.2 and theorem 2.7.11
in \cite{vw96}. (\ref{disine}) in T4 relates the Kullback-Leibler
divergence and Hellinger distance. Its proof depends on the
inequality that $\log x\leq 2(\sqrt{x}-1)$ for every $x>0$. T5 is
lemma 9.2 in \cite{mv01}. T6 is a result presented on page 79
of \cite{vg00} and is a special case of lemma 5.13 on the same
page, the proof of which can be found in pages 79--80. T7 and T8
are separately lemma 15.2 and 9.24 in \cite{kta}.
\end{remark}

{\it Proof of theorem~\ref{ratethm1}:} The definition of
$\hat{\eta}_{\tilde{\theta}_{n},\lambda_{n}}$ implies that
\begin{eqnarray*}
\lambda_{n}^{2}J^{2}(\hat{\eta}_{\tilde{\theta}_{n},\lambda_{n}})&\leq&\lambda_{n}^{2}J^{2}(\eta_{0})+(\mathbb{P}_{n}-P)
\left(\ell_{\tilde{\theta}_{n},\hat{\eta}_{\tilde{\theta}_{n},\lambda_{n}}}-\ell_{\tilde{\theta}_{n},\eta_{0}}\right)\\
&+&
P\left(\ell_{\tilde{\theta}_{n},\hat{\eta}_{\tilde{\theta}_{n},\lambda_{n}}}-\ell_{\tilde{\theta}_{n},\eta_{0}}\right)\\
&\leq&\lambda_{n}^{2}J^{2}(\eta_{0})+I+II.
\end{eqnarray*}
Note that by T6 and assumption (\ref{entrocon}), we have
\begin{eqnarray*}
I&\leq&(1+J(\hat{\eta}
_{\tilde{\theta}_{n},\lambda_{n}}))O_{P}(n^{-1/2})\times\left\{\left\|\frac{\ell_{\tilde{\theta}_{n},\hat{\eta}
_{\tilde{\theta}_{n},\lambda_{n}}}-\ell_{0}}{1+J(\hat{\eta}_{\tilde{\theta}_{n},\lambda_{n}})}\right\|_{2}^{1-\frac{1}{2k}}\vee
n^{-\frac{2k-1}{2(2k+1)}}
\right\}\\
&&+(1+J(\eta_{0}))O_{P}(n^{-1/2})\times\left\{\left\|\frac{\ell_{\tilde{\theta}_{n},\eta_{0}}
-\ell_{0}}{1+J(\eta_{0})}\right\|_{2}^{1-\frac{1}{2k}}\vee
n^{-\frac{2k-1}{2(2k+1)}} \right\}.
\end{eqnarray*}
By assumption (\ref{egsm1}), we have
\begin{eqnarray*}
II\aplt-d^{2}_{\tilde{\theta}_{n}}(\hat{\eta}_{\tilde{\theta}_{n},\lambda_{n}},\eta_{0})+\|\tilde{\theta}_{n}-\theta_{0}\|^{2}.
\end{eqnarray*}
Combining with the above, we can deduce that
\begin{eqnarray}
\hat{d}^{2}_{n}+\lambda_{n}^{2}\hat{J}_{n}^{2}
&\aplt&(1+\hat{J}_{n})O_{P}(n^{-1/2})
\times\left\{\left(\frac{\hat{d}_{n}+\|\tilde{\theta}_{n}-
\theta_{0}\|}{1+\hat{J}_{n}}\right)^{1-\frac{1}{2k}}\vee
n^{-\frac{2k-1}{2(2k+1)}}\right\}\nonumber\\
&+&(1+J_{0})O_{P}(n^{-1/2})
\times\left\{\left(\frac{\|\tilde{\theta}_{n}-
\theta_{0}\|}{1+J_{0}}\right)^{1-\frac{1}{2k}}\vee
n^{-\frac{2k-1}{2(2k+1)}}\right\}\nonumber\\
&+&\lambda_{n}^{2}J^{2}_{0}+\|\tilde{\theta}_{n}-\theta_{0}\|^{2},\label{tech2}
\end{eqnarray}
where
$\hat{d}_{n}=d_{\tilde{\theta}_{n}}(\hat{\eta}_{\tilde{\theta}_{n},\lambda_{n}},\eta_{0})$,
$J(\eta_{0})=J_{0}$ and
$\hat{J}_{n}=J(\hat{\eta}_{\tilde{\theta}_{n},\lambda_{n}})$. The
above inequality follows from assumption (\ref{egsm3}). Combining
all of the above inequalities, we can deduce that
\begin{eqnarray}
&&u_{n}^{2}=O_{P}(1)+O_{P}(1)u_{n}^{1-\frac{1}{2k}}\label{ratee2},\\
&&v_{n}=v_{n}^{-1}O_{P}(\|\tilde{\theta}_{n}-\theta_{0}\|^{2})+u_{n}^{1-\frac{1}{2k}}O_{P}(\lambda_{n})+O_{P}(n^{-\frac{1}{2}}\lambda_{n}^{-1}
\|\tilde{\theta}_{n}-\theta_{0}\|^{1-\frac{1}{2k}}),\label{ratee1}
\end{eqnarray}
where $u_{n}=(\hat{d}_{n} +\|\tilde{\theta}_{n}-
\theta_{0}\|)/(\lambda_{n}+\lambda_{n}\hat{J}_{n})$ and
$v_{n}=\lambda_{n}\hat{J}_{n}+\lambda_{n}$. The equation
(\ref{ratee2}) implies that $u_{n}=O_{P}(1)$. Inserting
$u_{n}=O_{P}(1)$ into (\ref{ratee1}), we can know that
$v_{n}=O_{P}(\lambda_{n}+\|\tilde{\theta}_{n}-\theta_{0}\|)$,
which implies $u_{n}$ has the desired order. This completes the
whole proof. $\Box$

{\it Proof of corollary~\ref{eg1rate}:} Conditions
(\ref{eg1bou})--(\ref{egsm1}) can be verified easily in this
example based on the arguments in theorem~\ref{ratethm1}
because $\ddot{\ell}_{\theta,f}$ has finite second moment, and
$p_{\theta,f}$ is bounded away from zero and infinity uniformly for
$(\theta,f)$ ranging over the whole parameter space. Note that
$d_{\theta}(f,f_{0})=\|p_{\theta,f}-p_{0}\|_{2}\apgt\|q_{\theta,f}-q_{\theta_{0},f_{0}}\|_{2}$
by Taylor expansion. Then by the assumption that $EVar(U|V)$ is
positive definite, we know that
$\|q_{\tilde{\theta}_{n},\hat{f}_{\tilde{\theta}_{n},\lambda_{n}}}-q_{\theta_{0},f_{0}}\|_{2}=O_{P}(\lambda_{n}+\|\tilde{\theta}_{n}-\theta_{0}\|)$
implies
$\|\hat{f}_{\tilde{\theta}_{n},\lambda_{n}}-f_{0}\|_{2}=O_{P}(\lambda_{n}+\|\tilde{\theta}_{n}-\theta_{0}\|)$.
Thus we only need to show that the $\epsilon$-bracketing entropy
number of the function class $\mathcal{O}$ defined below is of
order $\epsilon^{-1/k}$ to complete the proof of
(\ref{eg1pratre})--(\ref{eg1j}):
\begin{eqnarray*}
\mathcal{O}\equiv\left\{\frac{\ell_{\theta,f}(X)}{1+J(f)}:\|\theta-\theta_{0}\|\leq
C_{1},\|f-f_{0}\|_{\infty}\leq C_{1}, J(f)<\infty\right\},
\end{eqnarray*}
for some constant $C_{1}$. Note that $\ell_{\theta,f}(X)/(1+J(f))$
can be rewritten as:
\begin{eqnarray}
\Delta
A^{-1}\log\Phi\left(\bar{q_{\theta,f}}A\right)+(1-\Delta)A^{-1}\log\left(1-\Phi\left(
\bar{q_{\theta,f}}A\right)\right),
\end{eqnarray}
where $A=1+J(f)$ and $\bar{q}_{\theta,f}\in\mathcal{O}_{1}$, where
\begin{eqnarray*}
\mathcal{O}_{1}&\equiv&\left\{\frac{q_{\theta,f}(X)}{1+J(f)}:\|\theta-\theta_{0}\|\leq
C_{1},\|f-f_{0}\|_{\infty}\leq C_{1}, J(f)<\infty\right\},
\end{eqnarray*}
and where we know
$H_{B}(\epsilon,\mathcal{O}_{1},L_{2}(P))\aplt\epsilon^{-1/k}$ by
T1.

We next calculate the $\epsilon$-bracketing entropy number with
$L_{2}$ norm for the class of functions $R_{1}\equiv\{k_{a}(t):
t\mapsto a^{-1}\log\Phi(at)\;\mbox{for $a\geq 1$ and
$t\in\mathbb{R}$}\}$. By some analysis we know that $k_{a}(t)$ is
strictly decreasing in $a$ for $t\in\mathbb{R}$, and
$\sup_{t\in\mathbb{R}}|k_{a}(t)-k_{b}(t)|\aplt|a-b|$ because
$|\partial/\partial a(k_{a}(t))|$ is bounded uniformly over
$t\in\mathbb{R}$. In addition, we know that $\sup_{a,b\geq A_{0},
t\in\mathbb{R}}|k_{a}(t)-k_{b}(t)|\aplt A_{0}^{-1}$ because the
function $u\mapsto u\log\Phi(u^{-1}t)$ has bounded derivative for
$0<u\leq 1$ uniformly over $t\in\mathbb{R}$. The above two
inequalities imply that the $\epsilon$-bracketing number with
uniform norm is of order $O(\epsilon^{-2})$ for
$a\in[1,\epsilon^{-1}]$ and is 1 for $a>\epsilon^{-1}$. Thus we
know $H_{B}(\epsilon, R_{1}, L_{2})=O(\log\epsilon^{-2})$. By
applying a similar analysis to $R_{2}\equiv\{k_{a}(t): t\mapsto
a^{-1}\log(1-\Phi(at))\;\mbox{for $a\geq 1$ and
$t\in\mathbb{R}$}\}$, we obtain that $H_{B}(\epsilon, R_{2},
L_{2})=O(\log\epsilon^{-2})$. Combining this with T7 and T8, we
deduce that $H_{B}(\epsilon,\mathcal{O},L_{2})\aplt\epsilon^{-1/k}$.
This completes the proof of (\ref{eg1pratre})--(\ref{eg1j}).

For the proof of (\ref{eg10ratre}), we apply arguments similar
to those used in the proof of theorem~\ref{ratethm1} but after setting
$\lambda_{n}$, $J_{0}$ and $\hat{J}_{n}$ to zero in (\ref{tech2}).
Then we obtain the following equality:
$\hat{d}_{n}^{2}=O_{P}(n^{-2k/(2k+1)})+\|\tilde{\theta}_{n}-\theta_{0}\|^{2}+O_{P}(n^{-1/2})\|\tilde{\theta}_{n}-\theta_{0}\|^{1-1/2k}+
O_{P}(n^{-1/2})(\|\tilde{\theta}_{n}-\theta_{0}\|+\hat{d}_{n})^{1-1/2k}$.
By treating $\|\tilde{\theta}_{n}-\theta_{0}\|\leq n^{-k/(2k+1)}$
and $\|\tilde{\theta}_{n}-\theta_{0}\|> n^{-k/(2k+1)}$ differently
in the above equality, we obtain (\ref{eg10ratre}).$\Box$

{\it Proof of corollary~\ref{eg2rate}:} Lemma 7.1 in
\cite{mv99} establishes that
\begin{eqnarray}
\left\|p_{\tilde{\theta}_{n},\hat{\eta}_{\tilde{\theta}_{n},\lambda_{n}}}-p_{\theta_{0},\eta_{0}}\right\|_{2}
+\lambda_{n}J(\hat{\eta}_{\tilde{\theta}_{n},\lambda_{n}})=O_{P}(\lambda_{n}
+\|\tilde{\theta}_{n}-\theta_{0}\|)\label{coro21}
\end{eqnarray}
after choosing
\begin{eqnarray*}
m_{\theta,\lambda,\eta}=\log\frac{p_{\theta,\eta}+p_{\theta,\eta_{0}}}{2p_{\theta,\eta_{0}}}-\frac{1}{2}\lambda^{2}(J^{2}(\eta)-
J^{2}(\eta_{0}))
\end{eqnarray*}
in theorem~\ref{nuirate}. Note that the map $\theta\mapsto
p_{\theta,\eta_{0}}/f^{W,Z}(w,z)$ is uniformly bounded away from
zero at $\theta=\theta_{0}$ and continuous around a neighborhood
of $\theta_{0}$. Hence $m_{\theta,\lambda,\eta}$ is well defined.
Moreover,
$\mathbb{P}_{n}m_{\theta,\lambda,\hat{\eta}_{\theta,\lambda}}\geq
\mathbb{P}_{n}m_{\theta,\lambda,\eta_{0}}$ by the inequality that
$((p_{\theta,\eta}+p_{\theta,\eta_{0}})/2p_{\theta,\eta_{0}})^{2}\geq(p_{\theta,\eta}/p_{\theta,\eta_{0}})$.
(\ref{coro21}) now directly implies (\ref{eg2j}). For the proof of
(\ref{eg2pratre}), we need to consider the conclusion of lemma 7.4
(i), which states that
\begin{eqnarray}
\|p_{\theta,\eta}-p_{\theta_{0},\eta_{0}}\|_{2}\apgt(\|\theta-\theta_{0}\|\wedge
1+\||\eta-\eta_{0}|\wedge 1\|_{2})\wedge 1.\label{coro22}
\end{eqnarray}
Thus we have proved (\ref{eg2pratre}). For (\ref{eg1ratre}), we
just replace the $m_{\theta,\lambda,\eta}$ with
$m_{\theta,0,\eta}$ in the proof of lemma 7.1 in \cite{mv99}. Thus
we can show that
$d_{\theta}(\eta,\eta_{0})=\|p_{\theta,\eta}-p_{\theta_{0},\eta_{0}}\|_{2}$.
By combining lemma~\ref{lerate} and (\ref{coro22}), we know that
$\|\hat{\eta}_{\tilde{\theta}_{n}}-\eta_{0}\|_{2}=O_{P}(\delta_{n}+\|\tilde{\theta}_{n}-\theta_{0}\|)$,
for $\delta_{n}$ satisfying
$K(\delta_{n},\mathcal{S}_{\delta_{n}},L_{2}(P))\leq\sqrt{n}\delta_{n}^{2}$.
Note that $K(\delta,\mathcal{S}_{\delta},L_{2}(P))$ is as defined in
(\ref{thm2k}). By similar analysis as used in the proof of lemma 7.1
in \cite{mv99} and the strengthened assumption on $\eta$, we then
find that $K(\delta_{n},\mathcal{S}_{\delta_{n}},L_{2}(P))\aplt
\delta^{1-1/2k}_{n}$, which leads to the desired convergence rate
given in (\ref{eg1ratre}). $\Box$

{\it Proof of theorem~\ref{ratethm3}.} Note that
\begin{eqnarray*}
Pl^{\lambda_{n}}_{\tilde{\theta}_{n},\hat{\eta}_{\tilde{\theta}_{n},\lambda_{n}},h}&-&Pl^{\lambda_{n}}_{\theta_{0},\eta_{0},h}\\&=&
V(\tilde{\theta}_{n},\hat{\eta}_{\tilde{\theta}_{n},\lambda_{n}})h-
\mathbb{P}_{n}l^{\lambda_{n}}_{\tilde{\theta}_{n},\hat{\eta}_{\tilde{\theta}_{n},\lambda_{n}},h}-2\lambda_{n}^{2}\int
h^{(k)}(\hat{\eta}_{\tilde{\theta}_{n},\lambda_{n}}^{(k)}-\eta_{0}^{(k)})dz\\
&=&-(\mathbb{V}_{n}-V)(\tilde{\theta}_{n},\hat{\eta}_{\tilde{\theta}_{n},\lambda_{n}})h+2\lambda_{n}^{2}\int
h^{(k)}\eta_{0}^{(k)}dz\\
&=&-(\mathbb{V}_{n}-V)(\theta_{0},\eta_{0})h+o_{P}^{\ast}(n^{-1/2})+2\lambda_{n}^{2}\int
h^{(k)}\eta_{0}^{(k)}dz\\
&=&O_{P}(n^{-1/2})+2\lambda_{n}^{2}\int h^{(k)}\eta_{0}^{(k)}dz.
\end{eqnarray*}
The last two equalities in the above follow from assumptions
(\ref{con1th3}) and (\ref{con2th3}). The Fr\'{e}chet
differentiability of $V(\cdot,\cdot)$ at $(\theta_{0},\eta_{0})$
establishes that
\begin{eqnarray*}
Pl^{\lambda_{n}}_{\tilde{\theta}_{n},\hat{\eta}_{\tilde{\theta}_{n},\lambda_{n}},h}&-&Pl^{\lambda_{n}}_{\theta_{0},\eta_{0},h}
\\&=&\dot{V}(\tilde{\theta}_{n}-\theta_{0},
\hat{\eta}_{\tilde{\theta}_{n},\lambda_{n}}-\eta_{0})+o_{P}^{\ast}(\|\tilde{\theta}_{n}-\theta_{0}\|+d_{\tilde{\theta}_{n}}(\hat{\eta}_{\tilde
{\theta}_{n},\lambda_{n}},\eta_{0}) )\\&&-2\lambda_{n}^{2}\int
h^{(k)}(\hat{\eta}^{(k)}_{\tilde{\theta}_{n},\lambda_{n}}-\eta_{0}^{(k)})dz.
\end{eqnarray*}
Combining the above two sets of equations, we have, by the
linearity of $\dot{V}(\cdot,\cdot)$, established that
\begin{eqnarray*}
\dot{V}(0,\hat{\eta}_{\tilde{\theta}_{n},\lambda_{n}})=O_{P}(n^{-1/2})+O_{P}(\|\tilde{\theta}_{n}-\theta_{0}\|)+
2\lambda_{n}^{2}\int_{\mathcal{Z}}h^{(k)}\hat{\eta}^{(k)}_{\tilde{\theta}_{n},\lambda_{n}}dz.
\end{eqnarray*}
Now by the invertibility of $\dot{V}(0,\cdot)$, we can deduce
that $d_{\tilde{\theta}_{n}}(\hat{\eta}_{\tilde
{\theta}_{n},\lambda_{n}},\eta_{0})=O_{P}(n^{-1/2}+\|\tilde{\theta}_{n}-\theta_{0}\|+\lambda_{n}^{2}J^{2}(
\hat{\eta}_{\tilde{\theta}_{n},\lambda_{n}}))$. $\Box$

{\it Proof of theorem~\ref{penexp}.} We first show (\ref{penmle}),
and then we need to state one lemma before proceeding to the proof
of (\ref{penexpf}). For the proof of (\ref{penmle}), note that
\begin{eqnarray*}
0=\pp_{n}\dot{\ell}(\hat{\theta}_{
\lambda_{n}},\hat{\theta}_{\lambda_{n}},\hat{\eta}_{\lambda_{n}})+2\lambda_{n}^{2}\int_{\mathcal{Z}}
\hat{\eta}_{\lambda_{n}}^{(k)}(z)h_{0}^{(k)}(z)dz.
\end{eqnarray*}
Combining the third order Taylor expansion of
$\hat{\theta}_{\lambda_{n}}\mapsto\pp_{n}\dot{\ell}(\hat{\theta}_{
\lambda_{n}},\theta,\eta)$ around $\theta_{0}$, where
$\theta=\hat{\theta}_{\lambda_{n}}$ and
$\eta=\hat{\eta}_{\lambda_{n}}$, with conditions (\ref{emno1}) and
(\ref{emno2}), the first term in the right-hand-side
of the above displayed equality
equals
$\mathbb{P}_{n}\tilde{\ell}_{0}-\tilde{I}_{0}(\hat{\theta}_{\lambda_{n}}-\theta_{0})+O_{P}(\lambda_{n}+\|\hat{\theta}
_{\lambda_{n}}-\theta_{0}\|)^{2}$. By the inequality
$2\lambda_{n}^{2}\int_{\mathcal{Z}}
\hat{\eta}_{\lambda_{n}}^{(k)}(z)h_{0}^{(k)}(z)dz\leq\lambda_{n}^{2}(J^{2}(\hat{\eta}_{\lambda_{n}})+J^{2}(h_{0}))$
and assumption (\ref{jrate}), the second term in the right-hand-side
of the above equality is equal to $O_{P}(\lambda_{n}+\|\hat{\theta}
_{\lambda_{n}}-\theta_{0}\|)^{2}$. Combining everything, we
obtain the following:
\begin{eqnarray}
\;\;\;\;\;\;\;\;\;\frac{1}{\sqrt{n}}\sum_{i=1}^{n}\tilde{I}_{0}^{-1}\tilde{\ell}_{0}(X_{i})=\sqrt{n}(\hat{\theta}_{\lambda_{n}}-\theta_{0})+
O_{P}(n^{1/2}(\lambda_{n}+\|\hat{\theta}_{\lambda_{n}}-\theta_{0}\|)^{2}).\label{mleint}
\end{eqnarray}
The right-hand-side of (\ref{mleint}) is of the order
$O_{P}(\sqrt{n}\lambda_{n}^{2}+\sqrt{n}w_{n}(1+w_{n}+\lambda_{n}))$,
where $w_{n}$ represents
$\|\hat{\theta}_{\lambda_{n}}-\theta_{0}\|$. However, its
left-hand-side is trivially $O_{P}(1)$. Considering the fact that
$\sqrt{n}\lambda_{n}^{2}=o_{P}(1)$, we can deduce that
$\hat{\theta}_{\lambda_{n}}-\theta_{0}=O_{P}(n^{-1/2})$. Inserting
this into the previous display completes the proof of
(\ref{penmle}).

We next prove (\ref{penexpf}). Note that
$\hat{\theta}_{\lambda_{n}}-\theta_{0}=O_{P}(n^{-1/2})$. Hence the
order of the remainder terms in (\ref{emno1}) and (\ref{emno2})
becomes
$O_{P}(\lambda_{n}+\|\tilde{\theta}_{n}-\hat{\theta}_{\lambda_{n}}\|)^{2}$
and
$O_{P}(\lambda_{n}+\|\tilde{\theta}_{n}-\hat{\theta}_{\lambda_{n}}\|)$,
respectively. Expression~(\ref{lnplexpps}) in lemma~\ref{le4} below
implies that
\begin{eqnarray}
log pl_{\lambda_{n}}(\hat{\theta}_{\lambda_{n}})&=&\log
pl_{\lambda_{n}}(\theta_{0})+n(\hat{\theta}_{\lambda_{n}}-\theta_{0})^{T}\mathbb{P}_{n}\tilde{\ell}_{0}
\label{difint}\\&-&\frac{n}{2}(\hat{\theta}_{\lambda_{n}}-\theta_{0})^{T}\tilde{I}_{0}(\hat{\theta}_{\lambda_{n}}-\theta_{0})+
O_P(n^{1/2}\lambda_{n}^{2}).\nonumber
\end{eqnarray}
The difference between (\ref{difint}) and (\ref{lnplexpps})
generates
\begin{eqnarray*}
\log pl_{\lambda_{n}}(\tilde{\theta}_{n})&=&\log
pl_{\lambda_{n}}(\hat{\theta}_{\lambda_{n}})+n(\tilde{\theta}_{n}-\hat{\theta}_{\lambda_{n}})^{T}\left(\mathbb{P}_{n}\tilde{\ell}_{0}
-\tilde{I}_{0}(\hat{\theta}_{\lambda_{n}}-\theta_{0})\right)\\
&-&\frac{n}{2}(\tilde{\theta}_{n}-\hat{\theta}_{\lambda_{n}})^{T}\tilde{I}_{0}(\tilde{\theta}_{n}-\hat{\theta}_{\lambda_{n}})
+O_{P}(g_{\lambda_{n}}(\|\tilde{\theta}_{n}-\hat{\theta}
_{\lambda_{n}}\|)).
\end{eqnarray*}
(\ref{penexpf}) is now immediately obtained after considering
(\ref{penmle}). $\Box$

{\it Proof of theorem~\ref{pdist}.} Suppose that
$F_{\lambda_{n}}(\cdot)$ is the penalized posterior profile
distribution of $\sqrt{n}\varrho_{n}$ with respect to the prior
$\rho(\theta)$, where the vector $\varrho_{n}$ is defined as
$\tilde{I}_{0}^{1/2}(\theta-\hat{\theta}_{n})$. The parameter set
for $\varrho_{n}$ is $\Xi_{n}$. $F_{\lambda_{n}}(\cdot)$ can be
expressed as:
\begin{eqnarray}
\;\;\;\;\;\;\;F_{\lambda_{n}}(\xi)=\frac{\int_{\varrho_{n}\in(-\infty,n^{-1/2}\xi]\cap\Xi_{n}}\rho(\hat{\theta}_{\lambda_{n}}
+\tilde{I}_{0}^{-\frac{1}{2}}\varrho_{n} )
\frac{pl_{\lambda_{n}}(\hat{\theta}_{\lambda_{n}}+\tilde{I}_{0}^{-\frac{1}{2}}\varrho_{n}
)}{pl_{\lambda_{n}}(\hat{\theta}_{\lambda_{n}})}d\varrho_{n}}
{\int_{\varrho_{n}\in\Xi_{n}}\rho(\hat{\theta}_{\lambda_{n}}+\tilde{I}_{0}^{-\frac{1}{2}}\varrho_{n}
)\frac{pl_{\lambda_{n}}(\hat{\theta}_{\lambda_{n}}+\tilde{I}_{0}^{-\frac{1}{2}}\varrho_{n}
)}{pl_{\lambda_{n}}(\hat{\theta}_{\lambda_{n}})}d\varrho_{n}}.\label{pdistab}
\end{eqnarray}

Note that $d\varrho_{n}$ in the above is the short notation for
$d\varrho_{n1}\times\ldots\times d\varrho_{nd}$. To prove
theorem~\ref{pdist}, we first partition the parameter set
$\Xi_{n}$ as $\{\Xi_{n}\cap\{\|\varrho_{n}\|_{2}>
r_{n}\}\}\cup\{\Xi_{n}\cap\{\|\varrho_{n}\|_{2}\leq r_{n}\}\}$. By
choosing the proper order of $r_{n}$, we find the posterior mass
in the first partition region is of arbitrarily small order, as
verified in
lemma~\ref{pdist}.1 immediately below,
and the mass inside the second partition region
can be approximated by a stochastic polynomial in powers of
$n^{-1/2}$ with error of order dependent on the smoothing
parameter, as verified in lemma~\ref{pdist}.2 below.
This basic technique applies
to both the denominator and the numerator, yielding the quotient
series, which gives the desired result.

{\it lemma~\ref{pdist}.1.} Choose $r_{n}=o(n^{-1/3})$ and
$\sqrt{n}r_{n}\rightarrow\infty$. Under the conditions of
theorem~\ref{pdist}, we have
\begin{eqnarray}
\int_{\|\varrho_{n}\|>
r_{n}}\rho(\hat{\theta}_{\lambda_{n}}+\tilde{I}_{0}^{-\frac{1}{2}}\varrho_{n}
)\frac{pl_{\lambda_{n}}(\hat{\theta}_{\lambda_{n}}+\tilde{I}_{0}^{-\frac{1}{2}}\varrho_{n}
)}
{pl_{\lambda_{n}}(\hat{\theta}_{\lambda_{n}})}d\varrho_{n}=O_P(n^{-M}),\label{postout}
\end{eqnarray}
for any positive number $M$.

{\it Proof:} Fix $r>0$. We then have
\begin{eqnarray*}
\lefteqn{\int_{\|\varrho_{n}\|>
r}\rho(\hat{\theta}_{\lambda_{n}}+\tilde{I}_{0}^{-\frac{1}{2}}\varrho_{n}
)\frac{pl_{\lambda_{n}}(\hat{\theta}_{\lambda_{n}}+\tilde{I}_{0}^{-\frac{1}{2}}\varrho_{n}
)}{pl_{\lambda_{n}}(\hat{\theta}_{\lambda_{n}})}d\varrho_{n}}&&\\
&\mbox{\hspace{0.5in}}\leq&I\{\Delta_{\lambda_{n}}^{r}<
-n^{-\frac{1}{2}}\}\exp(-\sqrt{n})\int_{\Theta}
\rho(\theta)d\theta+I\{\Delta_{\lambda_{n}}^{r}\geq
-n^{-\frac{1}{2}}\},
\end{eqnarray*}
where
$\Delta_{\lambda_{n}}^{r}=sup_{\|\varrho_{n}\|>r}\Delta_{\lambda_{n}}(\hat{\theta}_{\lambda_{n}}+\varrho_{n}\tilde{I}_{0}^{-1/2})$.
Then by lemma 3.2 in \cite{ck06}, $I\{\Delta_{\lambda_{n}}^{r}\geq
-n^{-\frac{1}{2}}\}=O_{P}(n^{-M})$ for any fixed $r>0$. This
implies that there exists a positive decreasing sequence
$r_{n}=o(n^{-1/3})$ with $\sqrt{n}r_{n}\rightarrow\infty$ such
that (\ref{postout}) holds. $\Box$

{\it lemma~\ref{pdist}.2.} Choose $r_{n}=o(n^{-1/3})$ and
$\sqrt{n}r_{n}\rightarrow\infty$. Under the conditions of
theorem~\ref{pdist}, we have
\begin{eqnarray}
\int_{\|\varrho_{n}\|\leq
r_{n}}\left|\frac{pl_{\lambda_{n}}(\hat{\theta}_{\lambda_{n}}+\tilde{I}_{0}^{-\frac{1}{2}}\varrho_{n}
)}{pl_{\lambda_{n}}(\hat{\theta})}
\rho(\hat{\theta}_{\lambda_{n}}+\tilde{I}_{0}^{-\frac{1}{2}}\varrho_{n}
)-\exp \left(-\frac{n}{2}\varrho_{n}^{T}\varrho_{n}\right)
\rho(\hat{\theta}_{\lambda_{n}})\right|\nonumber\\\times
d\varrho_{n}=O_P(\lambda_{n}^{2}).\label{inmass}
\end{eqnarray}

{\it Proof:}  The posterior mass over the region
$\|\varrho_{n}\|_{2}\leq r_{n}$ is bounded by
\begin{eqnarray*}
\int_{\|\varrho_{n}\|_{2}\leq
r_{n}}\left|\frac{pl_{\lambda_{n}}(\hat{\theta}_{\lambda_{n}}+\tilde{I}_{0}^{-\frac{1}{2}}\varrho_{n}
)}{pl_{\lambda_{n}}(\hat{\theta}_{\lambda_{n}})}\rho(\hat{\theta}_{\lambda_{n}})-\exp
\left(-\frac{n}{2}\varrho_{n}^{T}\varrho_{n}\right)
\rho(\hat{\theta}_{\lambda_{n}})\right|d\varrho_{n}\;\;(\ast)\mbox{\hspace{0.5in}}\\
+\int_{\|\varrho_{n}\|_{2}\leq
r_{n}}\left|\frac{pl_{\lambda_{n}}(\hat{\theta}_{\lambda_{n}}+\tilde{I}_{0}^{-\frac{1}{2}}\varrho_{n}
)}{pl_{\lambda_{n}}(\hat{\theta}_{\lambda_{n}})}
\rho(\hat{\theta}_{\lambda_{n}}+\tilde{I}_{0}^{-\frac{1}{2}}\varrho_{n}
)-\frac{pl_{\lambda_{n}}(\hat{\theta}_{\lambda_{n}}+\tilde{I}_{0}^{-\frac{1}{2}}\varrho_{n}
)}{pl_{\lambda_{n}}(\hat{\theta}_{\lambda_{n}})
}\rho(\hat{\theta}_{\lambda_{n}})\right|d\varrho_{n}.\;\;(\ast\ast)
\end{eqnarray*}
By (\ref{penexpf}), we obtain
\begin{eqnarray*}
(\ast)&=&\int_{\|\varrho_{n}\|_{2}\leq
r_{n}}\left[\rho(\hat{\theta}_{\lambda_{n}})\exp
\left(-\frac{n\varrho_{n}^{T}\varrho_{n}}{2}\right)\left|\exp(O_P(g_{\lambda_{n}}(\|\varrho_{n}\|)))-1\right|\right]d\varrho_{n}.
\end{eqnarray*}
Obviously the order of $(\ast)$ depends on that of
$|\exp(O_{P}(g_{\lambda_{n}}(\|\varrho_{n}\|)))-1|$ for
$\lambda_{n}$ satisfying (\ref{smooth}) and $\|\varrho_{n}\|\leq
r_{n}$. In order to analyze its order, we partition the set
$\{\lambda_{n}=o_{P}(n^{-1/4})\;\mbox{and}\;\lambda_{n}^{-1}=O_{P}(n^{k/(2k+1)})
\}$ with the set $\{\lambda_{n}=O_{P}(n^{-1/3})\}$, i.e.
$U_{n}=\{\lambda_{n}=o_{P}(n^{-1/4})\;\mbox{and}\;\lambda_{n}^{-1}=O_{P}(n^{k/(2k+1)})
\}\cap\{\lambda_{n}=O_{P}(n^{-1/3})\}$ and
$L_{n}=\{\lambda_{n}=o_{P}(n^{-1/4})\;\mbox{and}\;\lambda_{n}^{-1}=O_{P}(n^{k/(2k+1)})
\}\cap\{\lambda_{n}=O_{P}(n^{-1/3})\}^{C}$. For the set $U_{n}$,
we have
$|\exp(O_{P}(g_{\lambda_{n}}(\|\varrho_{n}\|)))-1|=g_{\lambda_{n}}(\|\varrho_{n}\|)\times
O_{P}(1)$. For the set $L_{n}$, we have
$O_{P}(g_{\lambda_{n}}(\|\varrho_{n}\|))=O_{P}(n\|\varrho_{n}\|\lambda_{n}^{2}+n^{1/2}
\lambda_{n}^{2})$. We can take
$r_{n}=n^{-1-\delta}\lambda_{n}^{-2}$ for some $\delta>0$ such
that $\sqrt{n}r_{n}\rightarrow\infty$ and $r_{n}=o(n^{-1/3})$.
Then
$|\exp(O_{P}(g_{\lambda_{n}}(\|\varrho_{n}\|)))-1|=(n\|\varrho_{n}\|\lambda_{n}^{2}+n^{1/2}\lambda_{n}^{2})\times
O_{P}(1)$. Combining with the above, we know that
$(\ast)=O_{P}(\lambda_{n}^{2})$. By similar analysis, we can also
show that $(\ast\ast)$ has the same order. This completes the
proof of lemma~\ref{pdist}.2. $\Box$

We next start the formal proof of theorem~\ref{pdist}. By
considering both lemma~\ref{pdist}.1 and lemma~\ref{pdist}.2,
we know the denominator of (\ref{pdistab}) equals
\begin{eqnarray*}
\int_{\{\|\varrho_{n}\|_{2}\leq
r_{n}\}\cap\Xi_{n}}\left[\exp\left(-\frac{n}{2}\varrho_{n}^{T}\varrho_{n}\right)
\rho(\hat{\theta}_{\lambda_{n}})\right]d\varrho_{n}+O_{P}(\lambda_{n}^{2}).
\end{eqnarray*}
The first term in the above display equals
\begin{eqnarray*}
n^{-1/2}\rho(\hat{\theta}_{\lambda_{n}})\int_{\{\|u_{n}\|_{2}\leq
\sqrt{n}r_{n}\}\cap\sqrt{n}\Xi_{n}}e^{-u_{n}^{T}u_{n}/2}du_{n}
&=&n^{-1/2}\rho(\hat{\theta}_{\lambda_{n}})
\int_{\mathbb{R}^{d}}e^{-u_{n}^{T}u_{n}/2}du_{n}\\&+&O(\lambda_{n}^{2}),
\end{eqnarray*}
where $u_n=\sqrt{n}\varrho_{n}$. The above equality follows from
the inequality that $\int_{x}^{\infty}e^{-y^{2}/2}dy\leq
x^{-1}e^{-x^{2}/2}$ for any $x>0$. Consolidating the above
analyses, we deduce that the denominator of (\ref{pdistab}) equals
$n^{-\frac{1}{2}}\rho(\hat{\theta}_{\lambda_{n}})
(2\pi)^{d/2}\label{bot}+O_P(\lambda_{n}^{2})$. The same analysis
also applies to the numerator, thus completing the whole proof.
$\Box$

{\it Proof of corollary~\ref{pmome}:} We only show (\ref{pmon1})
in what follows. (\ref{pmon2}) can be verified similarly.
Showing (\ref{pmon1}) is equivalent to establishing
$\tilde{E}^{\lambda_{n}}_{\theta|x}(\varrho_{n})=O_{P}(\lambda_{n}^{2})$. Note that $\tilde{E}^{\lambda_{n}}_{\theta|x}(\varrho_{n})$ can be
written as:
\begin{eqnarray*}
\tilde{E}^{\lambda_{n}}_{\theta|x}(\varrho_{n})=\frac{\int_{\varrho_{n}\in\Xi_{n}}
\varrho_{n}\rho(\hat{\theta}_{\lambda_{n}}+\tilde{I}_{0}^{-\frac{1}{2}}\varrho_{n}
)
\frac{pl_{\lambda_{n}}(\hat{\theta}_{\lambda_{n}}+\tilde{I}_{0}^{-\frac{1}{2}}\varrho_{n}
)}{pl_{\lambda_{n}}(\hat{\theta}_{\lambda_{n}})}
d\varrho_{n}}{\int_{\varrho_{n}\in\Xi_{n}}\rho(\hat{\theta}_{\lambda_{n}}+\tilde{I}_{0}^{-\frac{1}{2}}\varrho_{n}
)
\frac{pl_{\lambda_{n}}(\hat{\theta}_{\lambda_{n}}+\tilde{I}_{0}^{-\frac{1}{2}}\varrho_{n}
)}{pl_{\lambda_{n}}(\hat{\theta}_{\lambda_{n}})}d\varrho_{n}}.
\end{eqnarray*}
By analysis similar to that
applied in the proof of theorem~\ref{pdist},
we know the denominator in the above display is
$n^{-1/2}(2\pi)^{d/2}\rho(\hat{\theta}_{\lambda_{n}})+O_P(\lambda_{n}^{2})$
and the numerator is a random vector of order
$O_P(n^{-1/2}\lambda_{n}^{2})$. This yields the conclusion.$\Box$

{\it Proof of theorem~\ref{pquant}.} Note that (\ref{ppostexp2})
implies
$\kappa_{n\alpha}=\tilde{I}_{0}^{-1/2}z_{\alpha}+O_{P}(n^{1/2}\lambda_{n}^{2})$,
for any $\xi<\alpha<1-\xi$, where $\xi\in(0,\frac{1}{2})$. Note
also that the $\alpha$-th quantile of a $d$ dimensional standard normal distribution,
$z_{\alpha}$, is not unique if $d>1$. The classical
Edgeworth expansion implies that
$P(n^{-1/2}\sum_{i=1}^{n}\tilde{I}_{0}^{-1/2}\tilde{\ell}_{0}(X_{i})\leq
z_{\alpha}+a_{n}(\alpha))=\alpha$, where
$a_{n}(\alpha)=O(n^{-1/2})$, for $\xi<\alpha<1-\xi$. Note that
$a_{n}(\alpha)$ is uniquely determined for each fixed $z_{\alpha}$
since $\tilde{\ell}_{0}(X_{i})$ has at least one absolutely
continuous component. Let
$\hat{\kappa}_{n\alpha}=\tilde{I}_{0}^{-1/2}z_{\alpha}+
(\sqrt{n}(\hat{\theta}_{\lambda_{n}}-\theta_{0})-
n^{-1/2}\sum_{i=1}^{n}\tilde{I}_{0}^{-1}\tilde{\ell}_{0}(X_{i})
)+\tilde{I}_{0}^{-1/2}a_{n}(\alpha)$. Then
$P(\sqrt{n}(\hat{\theta}_{\lambda_{n}}-\theta_{0})\leq\hat{\kappa}_{n\alpha})=
\alpha$. Combining with (\ref{penmle}), we obtain
$\hat{\kappa}_{n\alpha}=\kappa_{n\alpha}+O_P(n^{1/2}\lambda_{n}^{2})$.
The uniqueness of $\hat{\kappa}_{n\alpha}$ up to order
$O_P(n^{1/2}\lambda_n^2)$ follows from that
of $a_{n}(\alpha)$ for each chosen $z_{\alpha}$.$\Box$

{\it Proof of lemma~\ref{eg2le0}.} Assumptions S1 and S2 are
verified in lemma 5 of \cite{ck062}. For the verifications of the
assumption E1, we first show the asymptotic equicontinuity
condition (\ref{smcon1}). Without loss of generality, we assume
that $\lambda_{n}$ is bounded below by a multiple of
$n^{-k/(2k+1)}$ and bounded above by $n^{-1/4}$ in view of
(\ref{smooth}). Thus
\begin{eqnarray*}
P\left(\frac{\dot{\ell}(\theta_{0},\theta_{0},
\hat{f}_{\tilde{\theta}_{n},\lambda_{n}})-\dot{\ell}_{0}}
{n^{\frac{1}{4k+2}}(\lambda_{n}+\|\tilde{\theta}_{n}-\theta_{0}\|)}\right)^{2}\aplt\frac{\|\hat{f}_{\tilde{\theta}_{n},\lambda_{n}}
-f_{0}\|_{2}^{2}}{n^{\frac{1}{2k+1}}(\lambda_{n}+\|\tilde{\theta
}_{n}-\theta_{0}\|)^{2}}=O_{P}\left(n^{-\frac{1}{2k+1}}
\right),\nonumber
\end{eqnarray*}
where (\ref{eg1pratre}) implies the equality in the above
expression.

By (\ref{eg1j}), we know that
$J(\hat{f}_{\tilde{\theta}_{n},\lambda_{n}})=O_{P}(1+\|\tilde{\theta}_{n}-\theta_{0}\|/\lambda_{n})$
and $\|\hat{f}_{\tilde{\theta}_{n},\lambda_{n}}\|_{\infty}$ is
bounded by some constant, since $f\in\mathcal{H}_{k}^{M}$. We then
define the set $\mathcal{Q}_{n}$ as follows:
\begin{eqnarray*}
\left\{\frac{\dot{\ell}(\theta_{0},\theta_{0},f)-\dot{\ell}_{0}}{n^{\frac{1}{4k+2}}(\lambda_{n}+\|\theta-\theta_{0}\|)}:
J(f)\leq C_{n}( 1+\frac{\|\theta-\theta_{0}\|}{\lambda_{n}}),
\|f\|_{\infty}\leq M,
\|\theta-\theta_{0}\|\leq\delta\right\}\\\cap\left\{g\in L_{2}(P):
Pg^{2}\leq C_{n}n^{-\frac{1}{2k+1}}\right\},
\end{eqnarray*}
for some $\delta>0$. Obviously the function
$n^{-1/(4k+2)}(\dot{\ell}(\theta_{0},\theta_{0},
\hat{f}_{\tilde{\theta}_{n},\lambda_{n}})-\dot{\ell}_{0})/
(\lambda_{n}+\|\tilde{\theta}_{n}-\theta_{0}\|))\in\mathcal{Q}_{n}$
on a set of probability arbitrarily close to one, as
$C_{n}\rightarrow\infty$. If we can show
$\lim_{n\rightarrow\infty}E^{\ast}\|\mathbb{G}_{n}\|_{\mathcal{Q}_{n}}<\infty$
by T2, then assumption (\ref{smcon1}) is verified. Note that
$\dot{\ell}(\theta_{0},\theta_{0},f)$ depends on $f$ in a Lipschitz
manner. Consequently we can bound
$H_{B}(\epsilon,\mathcal{Q}_{n},L_{2}(P))$ by the product of some
constant and $H(\epsilon,\mathcal{R}_{n},L_{2}(P))$ in view of T3.
$\mathcal{R}_{n}$ is defined as
\begin{eqnarray*}
\{H_{n}(f):J(H_{n}(f))\aplt \lambda_{n}^{-1}n^{-1/(4k+2)},
\|H_{n}(f)\|_{\infty}\aplt \lambda_{n}^{-1}n^{-1/(4k+2)}\},
\end{eqnarray*}
where
$H_{n}(f)=f/(n^{1/(4k+2)}(\lambda_{n}+\|\theta-\theta_{0}\|))$.
By \cite{bs67}, we know that
\begin{eqnarray*}
H(\epsilon,\mathcal{R}_{n},L_{2}(P))\aplt(\lambda_{n}^{-1}n^{\frac{-1}{(4k+2)}})/\epsilon)^{1/k}.
\end{eqnarray*}
Note that $\delta_{n}=n^{-1/(4k+2)}$ and $M_{n}=n^{(2k-1)/(4k+2)}$
in T2. Thus by calculation we know that
$K(\delta_{n},\mathcal{Q}_{n},L_{2}(P))\aplt\lambda_{n}^{-1/2k}n^{-1/(4k+2)}$.
Then by T2 we can show that
$\lim_{n\rightarrow\infty}E^{\ast}\|\mathbb{G}_{n}\|_{\mathcal{Q}_{n}}<\infty$.

We next show (\ref{eg2em3}). It suffices to verify that
the sequence of classes of functions $\mathcal{V}_{n}$ is
$P$-Glivenko-Cantelli, where
$\mathcal{V}_{n}\equiv\{\ell^{(3)}(\bar{\theta}_{n},\tilde{\theta}_{n},\hat{f}_{\tilde{\theta}_{n},\lambda_{n}})(x)\}$,
for every random sequence $\bar{\theta}_{n}\rightarrow\theta_{0}$
and $\tilde{\theta}_{n}\rightarrow\theta_{0}$ in probability. A
Glivenko-Cantelli theorem for classes of functions that change
with $n$ is needed. By revising theorem 2.4.3 in \cite{vw96} with
minor notational changes, we obtain the following suitable
extension of the uniform entropy Glivenko-Cantelli theorem: Let
$\mathcal{F}_{n}$ be suitably measurable classes of functions with
uniformly integrable functions and
$H(\epsilon,\mathcal{F}_{n},L_{1}(\mathbb{P}_{n}))=o_{P}^{\ast}(n)$
for any $\epsilon>0$. Then
$\|\mathbb{P}_{n}-P\|_{\mathcal{F}_{n}}\rightarrow 0$ in
probability for every $\epsilon>0$. We then apply this revised
theorem to the set $\mathcal{F}_{n}$ of functions
$\ell^{(3)}(t,\theta,f)$ with $t$ and $\theta$ ranging over a
neighborhood of $\theta_{0}$ and $\lambda_{n}J(f)$ bounded by a
constant. By the form of $\ell^{(3)}(t,\theta,f)$, the entropy
number for $\mathcal{V}_{n}$ is equal to that of
\begin{eqnarray*}
\tilde{\mathcal{F}}_{n}\equiv
\{\phi(q_{t,f_{t}(\theta,f)}(x))R(q_{t,f_{t}(\theta,f)}(x)):(t,\theta)\in
V_{\theta_{0}}, \lambda_{n}J(f)\leq C, \|f\|_{\infty}\leq M\}.
\end{eqnarray*}
By arguments similar to those used in lemma 7.2 of \cite{mv99},
we know that
$\sup_{Q}H(\epsilon,\tilde{\mathcal{F}}_{n},
L_{1}(Q))\aplt(1+\lambda_{n}^{-1}/\epsilon)^{1/k}=o_{P}(n)$.
Moreover, the $\tilde{\mathcal{F}}_{n}$ are uniformly bounded since
$f\in\mathcal{H}_{k}^{M}$. Considering the fact that the
probability that
$\mathcal{V}_{n}$ is contained in $\tilde{\mathcal{F}}_{n}$ tends to
1, we have completed the proof of (\ref{eg2em3}).

For the proof of (\ref{eg2em1}), we only need to show that
$\mathbb{G}_{n}(\ddot{\ell}(\theta_{0},\tilde{\theta}_{n},\hat{f}_{\tilde{\theta}_{n},\lambda_{n}})-\ddot{\ell}_{0})=o_{P}(1)$
since $\ddot{\ell}_{0}(x)$ is uniformly bounded in $x$. Note that
we only need to show (\ref{eg2em1}) holds for
$\tilde{\theta}_{n}=\hat{\theta}_{n}+o(n^{-1/3})$ based on the
arguments in lemma~\ref{pdist}.2. We next show that
$\mathbb{G}_{n}(\ddot{\ell}(\theta_{0},\tilde{\theta}_{n},\hat{f}_{\tilde{\theta}_{n},\lambda_{n}})-\ddot{\ell}_{0})=o_{P}(1+n^{1/3}
\|\tilde{\theta}_{n}-\theta_{0}\|)=o_{P}(1)$. By the rate
assumptions R1, we have
\begin{eqnarray*}
P\left(\frac{\ddot{\ell}(\theta_{0},\tilde{\theta}_{n},\hat{f}_{\tilde{\theta}_{n},\lambda_{n}})-\ddot{\ell}_{0}}
{1+n^{1/3}\|\tilde{\theta}_{n}-\theta_{0}\|}\right)^{2}\aplt\frac{\|\tilde{\theta}_{n}-\theta_{0}\|^{2}
+\|\hat{f}_{\tilde{\theta}_{n},\lambda_{n}}-f_{0}\|_{2}^{2}}{(1+n^{1/3}\|\tilde{\theta}_{n}-\theta_{0}\|)^{2}}=O_{P}(n^{-1/2}).
\end{eqnarray*}
We next define $\bar{\cal Q}_{n}$ as follows:
\begin{eqnarray*}
\left\{\frac{\ddot{\ell}(\theta_{0},\theta,f)-\ddot{\ell}_{0}}{1+n^{1/3}\|\theta-\theta_{0}\|}:
J(f)\leq C_{n}( 1+\frac{\|\theta-\theta_{0}\|}{\lambda_{n}}),
\|f\|_{\infty}\leq M,
\|\theta-\theta_{0}\|<\delta\right\}\\\cap\left\{g\in L_{2}(P):
Pg^{2}\leq C_{n}n^{-\frac{1}{2}}\right\}.
\end{eqnarray*}
Obviously the function
$(\ddot{\ell}(\theta_{0},\tilde{\theta}_{n},\hat{f}_{\tilde{\theta}_{n},\lambda_{n}})-\ddot{\ell}_{0})/
(1+n^{1/3}\|\tilde{\theta}_{n}-\theta_{0}\|)\in\bar{\mathcal{Q}}_{n}$
on a set of probability arbitrarily close to one, as
$C_{n}\rightarrow\infty$. If we can show
$\lim_{n\rightarrow\infty}E^{\ast}\|\mathbb{G}_{n}\|_{\bar{\mathcal{Q}}_{n}}\rightarrow
0$ by T2, then the proof of (\ref{eg2em1}) is completed. Accordingly,
note that
$\ddot{\ell}(\theta_{0},\theta,f)$ depends on $(\theta,f)$ in a
Lipschitz manner. Consequently we can bound
$H_{B}(\epsilon,\bar{\mathcal{Q}}_{n},L_{2}(P))$ by the product of
some constant and
$(H(\epsilon,\bar{\mathcal{R}}_{n},L_{2}(P))+\log(1/\epsilon))$ in
view of T3. $\bar{\mathcal{R}}_{n}$ is defined as
\begin{eqnarray*}
\{H_{n}(f):J(H_{n}(f))\aplt 1+(n^{1/3}\lambda_{n})^{-1},
\|H_{n}(f)\|_{\infty}\aplt 1+(n^{1/3}\lambda_{n})^{-1}\},
\end{eqnarray*}
where $H_{n}(f)=f/(1+n^{1/3}\|\theta-\theta_{0}\|)$. By
\cite{bs67}, we know that
\begin{eqnarray*}
H(\epsilon,\bar{\mathcal{R}}_{n},L_{2}(P))\aplt((1+n^{-1/3}\lambda_{n}^{-1})/\epsilon)^{1/k}.
\end{eqnarray*}
Then by analysis similar to that used
in the proof of (\ref{smcon1}), we can show that
$\lim_{n\rightarrow\infty}E^{\ast}\|\mathbb{G}_{n}\|_{\bar{\mathcal{Q}}_{n}}\rightarrow
0$ in view of T2. This completes the proof of (\ref{eg2em1}).

For the proof of (\ref{eg2em2}), it suffices to show that
$\mathbb{G}_{n}(\ell_{t,\theta}(\theta_{0},\bar{\theta}_{n},\hat{f}_{\tilde{\theta}_{n},\lambda_{n}})-
\ell_{t,\theta}(\theta_{0},\theta_{0},f_{0}))=o_{P}(1)$ for
$\tilde{\theta}_{n}=\hat{\theta}_{n}+o(n^{-1/3})$ and for
$\bar{\theta}_{n}$ between $\tilde{\theta}_{n}$ and $\theta_{0}$,
in view of lemma~\ref{pdist}.2. Then we can show that
$\mathbb{G}_{n}(\ell_{t,\theta}(\theta_{0},\bar{\theta}_{n},\hat{f}_{\tilde{\theta}_{n},\lambda_{n}})-
\ell_{t,\theta}(\theta_{0},\theta_{0},f_{0}))=o_{P}(1+n^{1/3}\|\tilde{\theta}_{n}-\theta_{0}\|)=o_{P}(1)$
by similar analysis as used in the proof of (\ref{eg2em1}).$\Box$

{\it Proof of lemma~\ref{eg2le}.} By the assumption that
$\Delta_{\lambda_{n}}(\tilde{\theta}_{n})=o_{P}(1)$, we have
$\Delta_{\lambda_{n}}(\tilde{\theta}_{n})-\Delta_{\lambda_{n}}(\theta_{0})\geq
o_{P}(1)$. Thus the following inequality holds:
\begin{eqnarray*}
n^{-1}\sum_{i=1}^{n}\log\left[\frac{lik(\tilde{\theta}_{n},\hat{f}_{\tilde{\theta}_{n},\lambda_{n}},X_{i})}
{lik(\theta_{0},\hat{f}_{\theta_{0},\lambda_{n}},X_{i})}\right]
-n^{-1}\lambda_{n}^{2}[J^{2}(\hat{f}_{\tilde{\theta}_{n},\lambda_{n}})
-J^{2}(\hat{f}_{\theta_{0},\lambda_{n}})]\geq o_{P}(1)
\end{eqnarray*}
By considering assumption (\ref{jrate}), the above inequality
simplifies to
\begin{eqnarray*}
n^{-1}\sum_{i=1}^{n}\log\left[\frac{H(\tilde{\theta}_{n},\hat{f}_{\tilde{\theta}_{n},\lambda_{n}};X_{i})}
{H(\theta_{0},\hat{f}_{\theta_{0},\lambda_{n}};X_{i})}\right]\geq
o_{P}(1),
\end{eqnarray*}
where $H(\theta,f;X)=\Delta\Phi(C-\theta
U-f(V))+(1-\Delta)(1-\Phi(C-\theta U-f(V)))$. By arguments
similar to those used in lemma~\ref{eg2le0}.2 and by T5, we know
$H(\tilde{\theta}_{n},\hat{f}_{\tilde{\theta}_{n},\lambda_{n}};X_{i})$
belongs to some $P$-Donsker class. Combining the above conclusion
and the inequality $\alpha\log x\leq\log(1+\alpha\{x-1\})$ for
some $\alpha\in(0,1)$ and any $x>0$, we can show that
\begin{eqnarray}
P\log\left[1+\alpha\left(\frac{H(\tilde{\theta}_{n},\hat{f}_{\tilde{\theta}_{n},\lambda_{n}};X_{i})}
{H(\theta_{0},\hat{f}_{\theta_{0},\lambda_{n}};X_{i})}-1\right)\right]\geq
o_{P}(1).\label{leineq}
\end{eqnarray}
The remainder of the proof
follows the proof of lemma 6 in \cite{ck062}.$\Box$

{\it Proof of lemma~\ref{eg1le0}.} The maps (\ref{mapcon}) are
uniformly bounded since $F(\cdot)$, $\dot{F}(\cdot)$ and
$\ddot{F}(\cdot)$ are all uniformly bounded in
$(-\infty,+\infty)$. This completes the verifications of $S1$.
Note that $(W,Z)$ are in $[0,1]^{2}$ and $h_0(\cdot)$ is
intrinsically bounded over $[0,1]$. Hence we can show that the
Fr\'{e}chet derivatives of
$\eta\mapsto\ddot{\ell}(\theta_{0},\theta_{0},\eta)$ and
$\eta\mapsto\ell_{t,\theta}(\theta_{0},\theta_{0},\eta)$ for any
$\eta\in\mathcal{H}_{k}$ are bounded operators, from which we
can deduce that
$|\ddot{\ell}(\theta_{0},\theta_{0},\eta)(X)-\ddot{\ell}_{0}(X)|$
is bounded by the product of some integrable function and
$|\eta-\eta_{0}|(Z)$. This ensures (\ref{smcon2}) and
(\ref{smcon3}). For (\ref{nobias3}),
$P\dot{\ell}(\theta_{0},\theta_{0},\eta)$ can be written as
$P(F(\theta_{0}w+\eta_{0})-F(\theta_{0}w+\eta(z)))(w-h_{0}(z))$
since $P\dot{\ell}_{0}=0$. Note that
$P(w-h_{0}(z))\dot{F}(\theta_{0}w+\eta_{0}(z))(\eta-\eta_{0})(z)=0$.
This implies that
$P\dot{\ell}(\theta_{0},\theta_{0},\eta)=P(F(\theta_{0}w+\eta_{0})-F(\theta_{0}w+\eta(z))+
\dot{F}(\theta_{0}w+\eta_{0}(z))(\eta-\eta_{0})(z))(w-h_{0}(z))$.
However, by the common Taylor expansion, we have
$|F(\theta_{0}w+\eta)-F(\theta_{0}w+\eta_{0})-\dot{F}(\theta_{0}w+\eta_{0})(\eta-\eta_{0})|\leq\|\ddot{F}\|_{\infty}|\eta
-\eta_{0}|^{2}$. This proves (\ref{nobias3}).

We next verify assumption E1. For the asymptotic
equicontinuity condition (\ref{smcon1}), we first apply analysis
similar to that
used in the proof of lemma~\ref{eg2le0} to obtain
\begin{eqnarray*}
P\left(\frac{\dot{\ell}(\theta_{0},\theta_{0},\hat{\eta}_{\tilde{\theta}_{n},\lambda_{n}})-\dot{\ell}_{0}}
{n^{\frac{1}{4k+2}}(\lambda_{n}+\|\tilde{\theta}_{n}-\theta_{0}\|)}\right)^{2}\aplt
O_{P}\left(n^{-\frac{1}{2k+1}}\right).
\end{eqnarray*}
By lemma 7.1 in \cite{mv99}, we know that
$J(\hat{\eta}_{\tilde{\theta}_{n},\lambda_{n}})=O_{P}(1+\|\tilde{\theta}_{n}-\theta_{0}\|/\lambda_{n})$
and $\|\hat{\eta}_{\tilde{\theta}_{n},\lambda_{n}}\|_{\infty}$ is
bounded in probability by a multiple of
$J(\hat{\eta}_{\tilde{\theta}_{n},\lambda_{n}})+1$. Now we
construct the set $\tilde{\mathcal{Q}}_{n}$ as follows:
\begin{eqnarray*}
&&\left\{\frac{\dot{\ell}(\theta_{0},\theta_{0},\eta)-\dot{\ell}_{0}}{n^{\frac{1}{4k+2}}(\lambda_{n}+\|\theta-\theta_{0}\|)}:
J(\eta)\leq C_{n}( 1+\frac{\|\theta-\theta_{0}\|}{\lambda_{n}}),
\|\eta\|_{\infty}\leq C_{n}(1+J(\eta)),\right.\\&&\left.
\|\theta-\theta_{0}\|<\delta\right\}\cap\left\{g\in L_{2}(P):
Pg^{2}\leq C_{n}n^{-\frac{1}{2k+1}}\right\}.
\end{eqnarray*}
Clearly, the probability that the function
$n^{-1/(4k+2)}(\dot{\ell}(\theta_{0},\theta_{0},\hat{\eta}_{\tilde{\theta}_{n},\lambda_{n}})-\dot{\ell}_{0})/
(\lambda_{n}+\|\tilde{\theta}_{n}-\theta_{0}\|))\in\tilde{\mathcal{Q}}_{n}$
approaches 1 as $C_{n}\rightarrow\infty$. We next show that
$\lim_{n\rightarrow\infty}E^{\ast}\|\mathbb{G}_{n}\|_{\tilde{\mathcal{Q}}_{n}}<\infty$ by T2. Note that
$\dot{\ell}(\theta_{0},\theta_{0},\eta)$ depends on $\eta$
in a Lipschitz manner. Consequently, we can bound
$H_{B}(\epsilon,\tilde{\mathcal{Q}}_{n},L_{2}(P))$ by the product
of some constant and $H(\epsilon,\mathcal{R}_{n},L_{2}(P))$ in
view of T3, where $\mathcal{R}_{n}$ is as defined in the proof of
lemma~\ref{eg2le0}. By similar calculations as those performed in
lemma~\ref{eg2le0}, we can obtain
$K(\delta_{n},\tilde{\mathcal{Q}}_{n},L_{2}(P))\aplt\lambda_{n}^{-1/2k}n^{-1/(4k+2)}$. Thus
$\lim_{n\rightarrow\infty}E^{\ast}\|\mathbb{G}_{n}\|_{\tilde{\mathcal{Q}}_{n}}<\infty$, and (\ref{smcon1}) follows.

Next we define
$\bar{\mathcal{V}}_{n}\equiv\{\ell^{(3)}(\bar{\theta}_{n},\tilde{\theta}_{n},\hat{\eta}_{\tilde{\theta}_{n},\lambda_{n}})(x)\}$.
Similar arguments as those used
in the proof of lemma~\ref{eg2le0} can be
directly applied to the verification of (\ref{eg2em3}) in this second
model. By the form of $\ell^{(3)}(t,\theta,\eta)$, the
entropy number for $\bar{\mathcal{V}}_{n}$ is bounded above by
that of
$\bar{\mathcal{F}_{n}}\equiv\{\ddot{F}(tw+\eta(z)+(\theta-t)h_{0}(z)):(t,\theta)\in
V_{\theta_{0}}, \lambda_{n}J(\eta)\leq C_{n},
\|\eta\|_{\infty}\leq C_{n}(1+J(\eta))\}$. Similarly, we know
$\sup_{Q}H(\epsilon,\bar{\mathcal{V}}_{n},L_{1}(Q))\leq\sup_{Q}H(\epsilon,\bar{\mathcal{F}}_{n},L_{1}(Q))\aplt((1+\lambda_{n}^{-1})/\epsilon)^{1/k}=o_{P}(n)$.
Moreover, the $\bar{\mathcal{F}}_{n}$ are uniformly bounded. This
completes the proof for (\ref{eg2em3}).

The proof of (\ref{eg2em1}) and (\ref{eg2em2}) follows arguments
quite similar to those used
in the proof of lemma~\ref{eg2le0}. In other words, we
can show that
$\mathbb{G}_{n}(\ddot{\ell}(\theta_{0},\tilde{\theta}_{n},\hat{\eta}_{\tilde{\theta}_{n},\lambda_{n}})-\ddot{\ell}_{0})=o_{P}(1+n^{1/3}
\|\tilde{\theta}_{n}-\theta_{0}\|)=o_{P}(1)$ and
$\mathbb{G}_{n}(\ell_{t,\theta}(\theta_{0},\tilde{\theta}_{n},\hat{\eta}_{\tilde{\theta}_{n},\lambda_{n}})-\ell_{t,\theta}(\theta_{0},\theta_{0},\eta_{0}))=o_{P}(1+n^{1/3}
\|\tilde{\theta}_{n}-\theta_{0}\|)$. This concludes the proof.$\Box$

{\it Proof of lemma~\ref{eg1le}:} The proof of lemma~\ref{eg1le}
is analogous to that of lemma~\ref{eg2le}.$\Box$

\begin{lemma}\label{le4}
Assuming the assumptions in theorem~\ref{penexp}, we have
\begin{eqnarray}
log pl_{\lambda_{n}}(\tilde{\theta}_{n})&=&\log
pl_{\lambda_{n}}(\theta_{0})+n(\tilde{\theta}_{n}-\theta_{0})^{T}\mathbb{P}_{n}\tilde{\ell}_{0}
\label{lnplexpps}\\&-&\frac{n}{2}(\tilde{\theta}_{n}-\theta_{0})^{T}\tilde{I}_{0}(\tilde{\theta}_{n}-\theta_{0})
+O_P(g_{\lambda_{n}}(\|\tilde{\theta}_{n}-\hat{\theta}_{\lambda_{n}}\|)),\nonumber
\end{eqnarray}
for any $\tilde{\theta}_{n}=\theta_{0}+o_{P}(1)$.
\end{lemma}
{\it Proof.} $n^{-1}(\log
pl_{\lambda_{n}}(\tilde{\theta}_{n})-\log
pl_{\lambda_{n}}(\theta_{0}))$ is bounded above and below by
\begin{eqnarray*}
\pp_{n}(\ell(\tilde{\theta}_{n},\tilde{\theta}_{n},\hat{\eta}_{\tilde{\theta}_{n},\lambda_{n}})-
\ell(\theta_{0},\tilde{\theta}_{n},\hat{\eta}_{\tilde{\theta}_{n},\lambda_{n}}))-\frac{1}{n}\lambda_{n}^{2}
(J^{2}(\hat{\eta}_{\tilde{\theta}_{n},\lambda_{n}})-J^{2}(\eta_{\theta_{0}}(\tilde{\theta}_{n},
\hat{\eta}_{\tilde{\theta}_{n},\lambda_{n}})))
\end{eqnarray*}
and
\begin{eqnarray*}
\pp_{n}(\ell(\tilde{\theta}_{n},\theta_{0},\hat{\eta}_{\theta_{0},\lambda_{n}})-
\ell(\theta_{0},\theta_{0},\hat{\eta}_{\theta_{0},\lambda_{n}}))-\frac{1}{n}\lambda_{n}^{2}
(J^{2}(\eta_{\tilde{\theta}_{n}}(\theta_{0},
\hat{\eta}_{\theta_{0},\lambda_{n}}))-J^{2}(\hat{\eta}_{\theta_{0},\lambda_{n}})),
\end{eqnarray*}
respectively. By the third order Taylor expansion of
$\tilde{\theta}_{n}\mapsto\mathbb{P}_{n}\ell(\tilde{\theta}_{n},\theta,\eta)$
around $\theta_{0}$, for $\theta=\tilde{\theta}_{n}$ and
$\eta=\hat{\eta}_{\tilde{\theta}_{n},\lambda_{n}}$, and the above
empirical no-bias conditions (\ref{emno1}) and (\ref{emno2}), we
can find that the order of the difference between
$\pp_{n}(\ell(\tilde{\theta}_{n},\tilde{\theta}_{n},\hat{\eta}_{\tilde{\theta}_{n},\lambda_{n}})-
\ell(\theta_{0},\tilde{\theta}_{n},\hat{\eta}_{\tilde{\theta}_{n},\lambda_{n}}))$
and
$(\tilde{\theta}_{n}-\theta_{0})^{T}\mathbb{P}_{n}\tilde{\ell}_{0}-(\tilde{\theta}_{n}-\theta_{0})^{T}(\tilde{I}_{0}/2)(\tilde{\theta}_{n}-\theta_{0})$
is $O_{P}(n^{-1}g_{\lambda_{n}}(\|\tilde{\theta}_{n}
-\hat{\theta}_{\lambda_{n}}\|))$. By the inequality
$J^{2}(\eta_{t}(\theta,\eta))\leq
2J^{2}(\eta)+2(\theta-t)^{2}J^{2}(h_{0})$, we know that
$\lambda_{n}^{2}
(J^{2}(\hat{\eta}_{\tilde{\theta}_{n},\lambda_{n}})-J^{2}(\eta_{\theta_{0}}(\tilde{\theta}_{n},
\hat{\eta}_{\tilde{\theta}_{n},\lambda_{n}})))=O_{P}(\|\tilde{\theta}_{n}-
\hat{\theta}_{\lambda_{n}}\|+ \lambda_{n})^{2}$ provided
assumptions (\ref{smooth}) and (\ref{jrate}) hold. Similar
analysis also applies to the lower bound. This proves
(\ref{lnplexpps}).$\Box$

{\bf Acknowledgments.} The authors thank Dr. Joseph Kadane for
several insightful discussions.

\end{document}